\documentclass[11pt]{article}

\usepackage{amsmath,amssymb}

\usepackage{graphics}
\usepackage{tikz}

\newcommand{\beq}[1]{\begin{equation}\label{#1}}
\newcommand{\enq}{\end{equation}}
\newcommand{\raf}[1]{(\ref{#1})}
\newcommand{\qed}{\hfill\hbox{\rule{6pt}{6pt}} \medskip}
\newcommand{\proof}{\noindent {\bf Proof}.\ }

\newcommand{\ou}{\overline{u}}
\newcommand{\ov}{\overline{v}}
\newcommand{\ow}{\overline{w}}

\newcommand{\oy}{\overline{y}}

\newcommand{\ST}{\;|\;}

\newtheorem{theorem}{Theorem}
\newtheorem{lemma}[theorem]{Lemma}

\newtheorem{corollary}[theorem]{Corollary}

\newtheorem{definition}[theorem]{Definition}

\newtheorem{remark}[theorem]{Remark}

\def\problemDefInstance{Instance}

\def\problemDefQuestion{Question}

\def\problemDef#1#2#3{
\begingroup
\offinterlineskip
\centerline{\vbox{
\halign{\vrule \hskip1em ## \hfil & ## \hfil \hskip1em\vrule\cr
\noalign{\hrule}
\multispan2\vrule height1.3em depth.8em \hfil{\sc #1}\hfil\vrule\cr
\noalign{\hrule}
\multispan2\vrule height.8em \hfil \vrule \cr
{\bf \problemDefInstance{} :} &
\bgroup
\multiply\hsize by 2
\divide\hsize by 3
\parindent0pt
\baselineskip1.2em
\lineskip1pt
\lineskiplimit0pt
\vtop{#2}
\egroup\cr
\multispan2\vrule height.8em \hfil \vrule \cr
{\bf \problemDefQuestion{} :} &
\bgroup
\multiply\hsize by 2
\divide\hsize by 3
\parindent0pt
\baselineskip1.2em
\lineskip1pt
\lineskiplimit0pt
\vtop{#3}
\egroup\cr
\multispan2\vrule height.4em \hfil \vrule \cr
\noalign{\hrule}
}
}}
\endgroup
\vskip.5em
} 

\usepackage{mathrsfs}

\title{Separable discrete functions: recognition and sufficient conditions}
\author{
Endre Boros\thanks{
MSIS Dep. of RBS and RUTCOR, Rutgers University,
100 Rockafeller Road, Piscataway, NJ 08854-8054, USA.
(endre.boros@rutgers.edu)}
\and
Ond\v{r}ej \v{C}epek\thanks{
Department of Theoretical Informatics and Mathematical Logic, Charles University,
Malostransk\'{e} n\'{a}m. 25, 11800 Praha 1, Czech Republic.
(ondrej.cepek@mff.cuni.cz)}
\and
Vladimir Gurvich\thanks{
MSIS Dep. of RBS and RUTCOR, Rutgers University,
100 Rockafeller Road, Piscataway, NJ 08854-8054, USA;
National Research University Higher School of Economics, Moscow Russia.
(vladimir.gurvich@rutgers.edu)}
}
\date{\today{}}

\begin{document}
\maketitle

\begin{abstract}
A discrete function of $n$ variables is a mapping $g : X_1 \times \ldots \times X_n \rightarrow A$, where $X_1, \ldots, X_n$, and $A$ are arbitrary finite sets. Function $g$ is called {\em separable} if there exist $n$ functions $g_i : X_i \rightarrow A$ for $i = 1, \ldots, n$, such that for every input $x_1, \ldots ,x_n$ the function $g(x_1, \ldots, x_n)$ takes one of the values $g_1(x_1), \ldots ,g_n(x_n)$. Given a discrete function $g$, it is an interesting problem to ask whether $g$ is separable or not. Although this seems to be a very basic problem concerning discrete functions, the complexity of recognition of separable discrete functions of $n$ variables is known only for $n=2$. In this paper we will show that a slightly more general recognition problem, when $g$ is not fully but only partially defined, is NP-complete for $n \geq 3$. We will then use this result to show that the recognition of fully defined separable discrete functions is NP-complete for $n \geq 4$.

The general recognition problem contains the above mentioned special case for $n=2$. This case is well-studied in the context of game theory, where (separable) discrete functions of $n$ variables are referred to as (assignable) $n$-person game forms. There is a known sufficient condition for assignability (separability) of two-person game forms (discrete functions of two variables) called (weak) total tightness of a game form. This property can be tested in polynomial time, and can be easily generalized both to higher dimension and to partially defined functions. We will prove in this paper that weak total tightness implies separability for (partially defined) discrete functions of $n$ variables for any $n$, thus generalizing the above result known for $n=2$. Our proof is constructive. Using a graph-based discrete algorithm we show how for a given weakly totally tight (partially defined) discrete function $g$ of $n$ variables one can construct separating functions $g_1, \ldots ,g_n$ in polynomial time with respect to the size of the input function.

\end{abstract}

\newpage

\section{Introduction}
A discrete function of $n$ variables is a mapping $g : X_1 \times \ldots \times X_n \rightarrow A$, where $X_1, \ldots, X_n$, and $A$ are arbitrary finite sets. Discrete function $g$ is called {\em separable} if there exist $n$ separating functions of one variable each, $g_i : X_i \rightarrow A$ for $i = 1, \ldots, n$, such that

\medskip

$g(x) = g_1(x_1) ~ or ~ \ldots ~ or ~ g(x) = g_n(x_n)$
for every input
$x = (x_1, \ldots, x_n) \in X_1 \times \cdots \times X_n$.

\medskip

To see some simple examples for separable and non-separable discrete functions consider the following discrete functions of two variables. For $n=2$ we can interpret a discrete function as an array, where the rows are indexed by $X_1$, the columns by $X_2$, and the entries of the array are the $g(x_1,x_2)\in A$ values for $x_1\in X_1$ and $x_2\in X_2$.

\begin{equation}
\label{examples}
\left[ {\begin{array}{cc}
a & b \\
c & d
\end{array} } \right]
\hspace{1cm}
\left[ {\begin{array}{ccc}
a & a & c \\
a & b & b \\
c & b & c
\end{array} } \right]
\hspace{1cm}
\left[ {\begin{array}{ccc}
a & b & c\\
b & c & a
\end{array} } \right]
\hspace{1cm}
\left[ {\begin{array}{cccc}
a & a & b & b \\
c & d & c &d
\end{array} } \right]
\hspace{1cm}
\left[ {\begin{array}{ccc}
a & a & b \\
a & a & c \\
b & c & b
\end{array} } \right]
\end{equation}

The first two examples are clearly separable (we leave to the reader the easy task of defining the two separating functions, i.e. of assigning the correct values to rows and columns) while the last three examples are not (which is easy to prove by picking an entry in the array, assigning this value first to the corresponding row and then to the corresponding column, and following in both cases the sequence of forced assignments to a contradiction).

The concept of separability can be naturally extended to partially defined discrete functions of $n$ variables by requiring the above property only for those inputs for which the function is defined. Note that the separating functions can be always assumed to be fully defined - indeed if a set of $n$ partially defined separating functions fulfills the above condition, then so does any completion of this set to fully defined functions. Several concepts of decomposability for partially defined Boolean functions (a special subclass of discrete functions) are surveyed, for example, in \cite{BGHIK95}.

Although discrete functions are a general concept found in many areas of discrete mathematics, our motivation and interest in studying them came from game theory, and so we will in the rest of this paper switch to the standard game theoretical terminology. We denote by $I = [n] = \{1, \ldots, n\}$ the set of {\em players}, by $X_i$ the set of {\em strategies} of player $i \in I$, and by $A$ the set of possible {\em outcomes}.
A discrete function $g : X_1 \times \ldots \times X_n \rightarrow A$ of $n$ variables is called an $n$-person {\em game form}. In other words, once all players chose a particular strategy, say $x_i\in X_i$ for $i\in I$, then $g(x_1,...,x_n)$ denotes the corresponding outcome of the game. The  vector of the chosen particular strategies, $x = (x_1, \ldots, x_n) \in X_1 \times \cdots \times X_n$ is called a  {\em strategy profile}. In case $g$ is separable, the separating function $g_i$  {\em assigns} an outcome to every strategy of player $i$, and so in this context the term {\em assignability} (or property $AS$ in short) is used instead of separability.

The problem studied in this paper can be formulated as follows: given an $n$-person (partially defined) game form (an $n$-variate discrete function) decide whether it is assignable (separable). We shall see in Section~\ref{complexity-section} that this problem is solvable in polynomial time for $n=2$ and that it is computationally hard for $n \geq 3$ for the partially defined case, and for $n \geq 4$ for the fully defined case. In those situations when recognizing assignability is difficult, it makes sense to look for conditions that (i) are sufficient for assignability and (ii) can be tested in polynomial time. Such conditions will be studied in Section~\ref{main-section}. First, let us recall results related to the case of two players. This case has been most studied and best understood so far.

\subsection{The two-person case}
Since assignability is obviously a hereditary property (any subform of an assignable game form, induced by subsets of the strategy sets of the players, is obviously also assignable), a natural question arises, whether assignable game forms admit a characterization by a finite set of forbidden subforms.
Interestingly, this is not the case already for two-person game forms, as the following infinite family (from \cite{BCGMZ09}) of non-assignable game forms demonstrates it:

\begin{equation}
\label{sequence}
\left[ {\begin{array}{ccc}
a & b & a \\
c & a & b \\
b & b & a
\end{array} } \right]
\hspace{1cm}
\left[ {\begin{array}{cccc}
a & b & a & a \\
c & a & b & a \\
c & b & a & b \\
b & b & b & a
\end{array} } \right]
\hspace{1cm}
\left[ {\begin{array}{ccccc}
a & b & a & a & a \\
c & a & b & a & a \\
c & b & a & b & a \\
c & b & b & a & b \\
b & b & b & b & a
\end{array} } \right]
\hspace{1cm}
\cdots
\end{equation}

It is not hard to see, that each such game form is not assignable, while a removal of any row or column makes it assignable, i.e. the presented game forms are minimal non-assignable. On the other hand, it can be shown \cite{BCGMZ09}, that if assignability is made more restrictive by requiring that every strategy profile is covered by exactly one of the players (either $g(x_1, x_2) = g_1(x_1)$ or $g(x_1, x_2) = g_2(x_2)$  but not both), then such {\em strong assignability} can be characterized by a finite set of minimal non-assignable game forms. Note also, that any two person game form with at most two outcomes is clearly assignable (assign one outcome to all rows and the other outcome to all columns), and hence the game forms in~(\ref{sequence}) use the minimal number of three outcomes needed for non-assignability.

This observation can be generalized to the case of $n$ players: any $n$-person game form with at most $n$ outcomes is clearly assignable, and it is easy to construct non-assignable $n$-person game forms with $n+1$ outcomes. Such a construction proceeds as follows. Take an $m \times m \times \cdot \times m$ game form, and for each of the first $n$ outcomes select $m$ strategy profiles (thus we need $nm < m^n$ to have enough strategy profiles) such that no two of them share any coordinate value. In the two player case these are two permutation submatrices, in~(\ref{sequence}) the main diagonal with outcome $a$ and the diagonal above it plus the bottom left corner with outcome $b$. Each of these $m$-tuples of strategy profiles requires $m$ different player strategies to cover all $m$ outcomes. For each outcome $k$ we need $g_i(x_j)=k$ for $m$ distinct pairs $(i,j)$), and so altogether they use all $nm$ available strategies $g_i(x_j)$. Hence, putting the last outcome $n+1$ to any still "vacant" strategy profile causes the game form to be non-assignable. Let us finally remark, that this construction works for any $n$ and $m$ such that $nm < m^n$, but to make such a game form minimal non-assignable is probably a very tricky task (we do not see any easy way how to generalize the construction of minimal non-assignable game forms to any $n>2$).

Several other properties of two-person game forms which are connected to assignability were studied in the literature. An important property is {\em tightness}. We say that a player can \emph{guarantee} a subset $B\subseteq A$ of outcomes, if he/she has a strategy such that no matter what strategy the other player chooses, the corresponding outcome belongs to $B$. A two-person game form $g : X_1 \times X_2 \rightarrow A$ is called {\em tight} (has property $T$, in short) if for any subset $B \subseteq A$ either one player can guarantee $B$ or the other player can guarantee $\bar{B}=A\setminus B$. Note that both of the above cannot happen, while there may be a subset $B\subseteq A$ in a non-tight game form such that neither player $1$ can guarantee $B$, nor player $2$ can guarantee $\bar{B}$. The importance of tightness stems partly from the fact that it is in the two player case equivalent to Nash-solvability of a game form \cite{Gur75,Gur78,Gur88}, which is a pivotal concept in game theory. No polynomial time algorithm for verifying tightness is known, however a quasi-polynomial one was suggested
by Fredman and Khachiyan in \cite{FK96}.

Another property of game forms related to tightness is total tightness.
A two-person game form $g : X_1 \times X_2 \rightarrow A$ is called {\em totally tight} ($TT$) if every $2 \times 2$ subform of $g$ (which is a two-dimensional array in this case) is tight, or equivalently, if it contains a constant line (i.e. row or column).
More precisely, let us call $g' : X'_1 \times X'_2 \rightarrow A$
to be a $2 \times 2$ restriction of $g$ if $X'_1 = \{x_1, x'_1\} \subseteq X_1$ and $X'_2 = \{x_2, x'_2\} \subseteq X_2$ are 2-element subsets of $X_1$ and $X_2$.
Then $g$ is $TT$ if for every $2 \times 2$ restriction $g'$ of $g$ we have
\[ g'(x_1, x_2) = g'(x_1, x'_2), \mbox{ or }
g'(x_1, x_2) = g'(x'_1,x_2), \mbox{ or }
g'(x'_1, x'_2) = g'(x'_1,x_2), \mbox{ or }
g'(x'_1, x'_2) = g'(x_1, x'_2). \]

It was shown in \cite{BGMP10} that totally tight two-person game forms are both tight and assignable.
It is easy to observe, that checking whether a two-person game form is $TT$ can be done in polynomial time, as there is only a polynomial number of $2 \times 2$ subforms to check (this easy fact was observed e.g. in~\cite{BCG11} ). Thus property $TT$ provides a simple sufficient condition for the assignability of two-person game forms.

\subsection{The $n$-person case for $n \geq 3$}
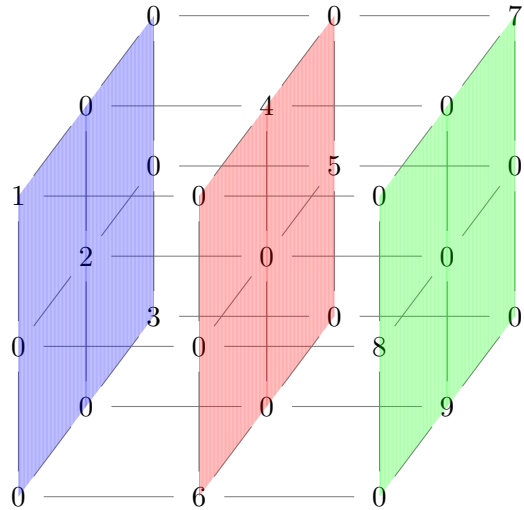
\begin{figure}[ht]
\centering
\begin{tikzpicture}[z={(0.45cm,0.6cm)},x={(1.2cm,0cm)},y={(0cm,1cm)},scale=2]

\foreach \x in {1,2,3} \foreach \y in {1,2,3}
	\draw[gray,thin] (\x,\y,1) -- (\x,\y,3);

\foreach \x in {1,2,3} \foreach \z in {1,2,3}
	\draw[gray,thin] (\x,1,\z) -- (\x,3,\z);
	
\foreach \y in {1,2,3} \foreach \z in {1,2,3}
	\draw[gray,thin] (1,\y,\z) -- (3,\y,\z);

\node[circle,fill=white] at (1,1,1) {$0$};
\node[circle,fill=white] at (1,2,1) {$0$};
\node[circle,fill=white] at (1,3,1) {$1$};
\node[circle,fill=white] at (1,1,2) {$0$};
\node[circle,fill=white] at (1,2,2) {$2$};
\node[circle,fill=white] at (1,3,2) {$0$};
\node[circle,fill=white] at (1,1,3) {$3$};
\node[circle,fill=white] at (1,2,3) {$0$};
\node[circle,fill=white] at (1,3,3) {$0$};

\node[circle,fill=white] at (2,1,1) {$6$};
\node[circle,fill=white] at (2,2,1) {$0$};
\node[circle,fill=white] at (2,3,1) {$0$};
\node[circle,fill=white] at (2,1,2) {$0$};
\node[circle,fill=white] at (2,2,2) {$0$};
\node[circle,fill=white] at (2,3,2) {$4$};
\node[circle,fill=white] at (2,1,3) {$0$};
\node[circle,fill=white] at (2,2,3) {$5$};
\node[circle,fill=white] at (2,3,3) {$0$};

\node[circle,fill=white] at (3,1,1) {$0$};
\node[circle,fill=white] at (3,2,1) {$8$};
\node[circle,fill=white] at (3,3,1) {$0$};
\node[circle,fill=white] at (3,1,2) {$9$};
\node[circle,fill=white] at (3,2,2) {$0$};
\node[circle,fill=white] at (3,3,2) {$0$};
\node[circle,fill=white] at (3,1,3) {$0$};
\node[circle,fill=white] at (3,2,3) {$0$};
\node[circle,fill=white] at (3,3,3) {$7$};

\shade[right color=blue,left color=blue,opacity=0.2] (1,1,1) -- (1,3,1) -- (1,3,3) -- (1,1,3) -- cycle;

\shade[right color=red,left color=red,opacity=0.2] (2,1,1) -- (2,3,1) -- (2,3,3) -- (2,1,3) -- cycle;

\shade[right color=green,left color=green,opacity=0.2] (3,1,1) -- (3,3,1) -- (3,3,3) -- (3,1,3) -- cycle;

\end{tikzpicture}
\caption{A non-assignable 3D example.\label{f-no-3D}}
\end{figure}

As we shall recall in Section \ref{complexity-section} that assignability of two-person game forms can be tested in polynomial time. A natural idea is to try to reduce the assignability of higher dimensional game forms to the two-dimensional case by considering projections. Let  $g : X_1 \times \ldots \times X_n \rightarrow A$ be a $n$-person game form. An $i$-th two dimensional projection $g^i$ of $g$ is a two person game form where the set of strategies of the first player consists of the strategies of player $i$ and the set of strategies of the second player is the direct product of all the strategies of the remaining $n-1$ players. Is it true that $g$ is assignable if and only if $g^i$ is assignable for every $i$? Unfortunately not. Both implications fail already for $n=3$. First we shall show an example of a $3 \times 3 \times 3$ three-person game form which is not assignable, but each of its three two-dimensional $3\times 9$  projections is assignable. The game form is given in Figure \ref{f-no-3D}.

\begin{figure}[htb]
\centering
\begin{tikzpicture}[z={(0.45cm,0.6cm)},x={(-1.2cm,0cm)},y={(0cm,1cm)},scale=2]

\foreach \x in {1,2,3} \foreach \y in {1,2,3}
	\draw[gray,thin] (\x,\y,1) -- (\x,\y,3);

\foreach \x in {1,2,3} \foreach \z in {1,2,3}
	\draw[gray,thin] (\x,1,\z) -- (\x,3,\z);
	
\foreach \y in {1,2,3} \foreach \z in {1,2,3}
	\draw[gray,thin] (1,\y,\z) -- (3,\y,\z);

\node[circle,fill=white] at (1,1,1) {$c$};
\node[circle,fill=white] at (1,2,1) {$b$};
\node[circle,fill=white] at (1,3,1) {$a$};
\node[circle,fill=white] at (1,1,2) {$a$};
\node[circle,fill=white] at (1,2,2) {$c$};
\node[circle,fill=white] at (1,3,2) {$b$};
\node[circle,fill=white] at (1,1,3) {$*$};
\node[circle,fill=white] at (1,2,3) {$a$};
\node[circle,fill=white] at (1,3,3) {$c$};

\node[circle,fill=white] at (2,1,1) {$*$};
\node[circle,fill=white] at (2,2,1) {$*$};
\node[circle,fill=white] at (2,3,1) {$b$};
\node[circle,fill=white] at (2,1,2) {$*$};
\node[circle,fill=white] at (2,2,2) {$*$};
\node[circle,fill=white] at (2,3,2) {$c$};
\node[circle,fill=white] at (2,1,3) {$*$};
\node[circle,fill=white] at (2,2,3) {$*$};
\node[circle,fill=white] at (2,3,3) {$*$};

\node[circle,fill=white] at (3,1,1) {$*$};
\node[circle,fill=white] at (3,2,1) {$*$};
\node[circle,fill=white] at (3,3,1) {$c$};
\node[circle,fill=white] at (3,1,2) {$*$};
\node[circle,fill=white] at (3,2,2) {$*$};
\node[circle,fill=white] at (3,3,2) {$a$};
\node[circle,fill=white] at (3,1,3) {$*$};
\node[circle,fill=white] at (3,2,3) {$*$};
\node[circle,fill=white] at (3,3,3) {$*$};

\shade[right color=blue,left color=blue,opacity=0.2] (1,1,1) -- (1,1,2) -- (1,3,2) -- (1,3,1) -- cycle;

\shade[right color=red,left color=red,opacity=0.2] (1,2,1) -- (1,2,3) -- (1,3,3) -- (1,3,1) -- cycle;

\shade[right color=green,left color=green,opacity=0.2] (1,3,1) -- (1,3,2) -- (3,3,2) -- (3,3,1) -- cycle;

\end{tikzpicture}
\caption{Assignable 3D example that has no 2D assignable projection, where symbol $*$ can be any of $a$, $b$ or $c$.\label{f-3D-no-2D}}
\end{figure}
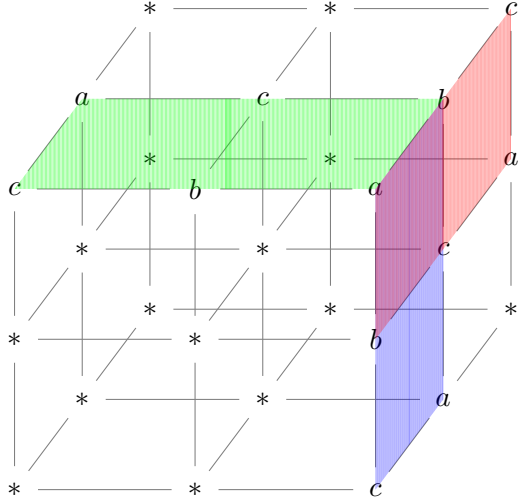

Clearly, this three person game form is not assignable as there are 10 outcomes but only 9 values to be assigned. On the other hand, each of the three two-dimensional projections is a $9 \times 3$ two-person game form which is assignable by assigning $0$ to all three columns and assigning outcomes $1$ to $9$ to the rows (note that each row contains exactly two $0$'s and exactly one other outcome in each of the three projections). Second we shall show an example of a $3 \times 3 \times 3$ three-person game form which is assignable, but none of its three two-dimensional projections is assignable. The game form is given in Figure \ref{f-3D-no-2D}.
In this example the symbol $*$ stands for an arbitrary outcome from the set $\{a,b,c\}$. Since there are only three outcomes, the three-person game form is trivially assignable by assigning outcome $a$ to planes orthogonal to direction $1$, outcome $b$ to planes orthogonal to direction $2$, and outcome $c$  to planes orthogonal to direction $3$. On the other hand, each of the three $9\times 3$ two-dimensional projections is a game form that contains a $2 \times 3$ subform
\[
\hspace{1cm}
\left[ {\begin{array}{ccc}
a & b & c\\
b & c & a
\end{array} } \right]
\hspace{1cm}
\]
which is not assignable, and hence none of the projections is assignable (note that this $2 \times 3$ subform is the middle example in~(\ref{examples})).

Therefore, to obtain sufficient conditions for assignability for $n \geq 3$ we need to look at other properties of game forms, e.g. generalizations of the properties studied in the two-person case.

The concepts of tightness and total tightness can be extended to the general $n$-person case in the following way.
Given an $n$-person game form $g : X_1 \times \ldots \times X_n \rightarrow A$ and a partition $I = K \cup \overline{K}$ of the players
into two complementary non-empty coalitions, the two-person game form $g^K : X_K \times X_{\overline{K}} \rightarrow A$ is defined as follows.
The strategies of the first and second players are the elements of the direct products $X_K = \prod_{i \in K} X_i$  and
$X_{\overline{K}} = \prod_{i \in \overline{K}} X_i$, respectively. For $x \in X_K$ and $\overline{x} \in X_{\overline{K}}$ we define $g^K(x,\overline{x}) = g(y)$
where $y$ originates from $x$ and $\overline{x}$ by concatenating them and reordering the coordinates according to the $X_1, \ldots, X_n$ order.
An $n$-person game form $g$ is called {\em tight} (respectively, {\em totally tight}) if $g^K$ is tight (respectively, {\em totally tight}) for all non-empty $K \subseteq I$.

Similarly, we call $g$ {\em weakly tight} (respectively, {\em weakly totally tight} if $g^K$ is tight (respectively, totally tight) for all $K \subseteq I$ such that $|K| = 1$.
Note, that in this case we consider exactly the two-dimensional projections defined above in the first paragraph of this subsection.
We shall denote these concepts by $WT$ and $WTT$, respectively. Let us remark that, by the above definition, $T = WT$ and $TT = WTT$ for $n \leq 3$.
Indeed, in this case for any non-trivial partition one of the two coalitions contains only one player. Let us also remark, that the above defined concepts ($T$, $WT$, $TT$, and $WTT$)
can be straightforwardly extended to partially defined game forms by replacing all undefined values of the game form by a single extra outcome, and by requiring
the corresponding property for the resulting fully defined game form.

Let us observe that testing whether a three-person game form $g$ is (weakly) totally tight can be done in polynomial time. In fact, as we shall show in Section~\ref{complexity-section}, property $WTT$ can be tested in polynomial time (with respect to the size of $g$) for any $n$.
Therefore, property $WTT$ may be a good candidate for a property that we are looking for, a property which can be tested in polynomial time and which implies assignability. This leads to another results of this paper, namely that weak total tightness implies assignability of $n$-person game forms (both partially and fully defined) for all $n$.

The structure of the paper is as follows.  In Section \ref{main-section} we prove that $WTT \Rightarrow AS$ for every $n$, that is, that property $WTT$ implies assignability for every $n$. In Section \ref{complexity-section} we prove that deciding assignability of partially defined game forms is NP-complete for $n\geq 3$, and that deciding assignability of fully defined game forms is NP-complete for $n\geq 4$. We close the paper by providing further connections to game theory in Section \ref{game-theory}, and by listing some open problems in Section \ref{conclusions}.

\bigskip

\section{Weak total tightness implies assignability}
\label{main-section}
When dealing with game forms it is sometimes convenient to think of the $n$-person game form as of an $n$-dimensional array and use a geometric interpretation for subarrays. A {\em line in direction $i$} is a set of strategy profiles (a $1$-dimensional subarray) where all coordinates are fixed and only coordinate $i$ is used as a running index. In game theoretic terms a line in direction $i$ is a $1$-dimensional subform obtained by fixing the strategies of all players except of player $i$. A {\em hyperplane perpendicular to direction $i$} is a set of strategy profiles (an $(n-1)$-dimensional subarray) where all coordinates are used as running indices and only coordinate $i$ is fixed. In game theoretic terms a hyperplane perpendicular to direction $i$ is an $(n-1)$-dimensional subform obtained by fixing the strategy of player $i$.

\begin{definition}\label{constant}
Given a game form $g$, a set $S \subseteq X$ of strategy profiles will be called a {\em constant region} if all strategy profiles in $S$ get the same outcome, i.e. if there exists an outcome $c \in A)$ such that for all strategy profiles $x \in S$ we have $g(x)=c$.
\end{definition}

\begin{remark}\label{basic-assumptions}
We will assume in the remainder of this paper that the game form we are dealing with contains no constant hyperplane (a hyperplane which is a constant region) and no pair of duplicate parallel hyperplanes. These assumptions can be made without a loss of generality as such hyperplanes obviously influence neither total tightness nor assignability of the considered game forms.
\end{remark}

Let us also note that for a $WTT$ game form $g$ the entries in any $2 \times 2$ subarray of $g_{\{i\}}$ (recall that $g_{\{i\}}$ is a two player game form where player $i$ constitutes a one person coalition and all other players constitute the complementary coalition) can be geometrically thought of as the four intersections of two arbitrary distinct lines in direction $i$ with two arbitrary distinct $(n-1)$-dimensional hyperplanes perpendicular to direction $i$.

Let us start with a simple lemma describing a forbidden substructure for WTT game forms.

\begin{lemma}
\label{forbidden}
Let $g$ be an $n$-person game form, $i$ an arbitrary direction (player), and $H_j, H_k$ be two distinct parallel hyperplanes perpendicular to direction $i$. Furthermore let $\ell^i, 1\leq i\leq 4$ be four lines (not necessarily all distinct) in direction $i$ intersecting $H_j$ in strategy profiles $x_j^i, 1\leq i\leq 4$ and intersecting $H_k$ in strategy profiles $x_k^i, 1\leq i\leq 4$, such that
\begin{description}
	\item[(1)] $g(x_j^1) \neq g(x_j^2)$,
	\item[(2)] $g(x_k^3) \neq g(x_k^4)$,
	\item[(3)] $g(x_j^i) \neq g(x_k^i), 1\leq i\leq 4$.
\end{description}
Then $g$ is not a WTT game form.
\end{lemma}

\begin{figure}[ht]
\centering
\begin{tikzpicture}[x={(4cm,0cm)},y={(0cm,2cm)}]

\draw[thin,step=1] (-0.5,-0.5) grid (1.5,3.5);

\node[draw,circle,fill=white] at (0,0) {\footnotesize $g(x^4_j)$};
\node[draw,circle,fill=white] at (0,1) {\footnotesize $g(x^3_j)$};
\node[draw,circle,fill=white] at (0,2) {\footnotesize $g(x^2_j)$};
\node[draw,circle,fill=white] at (0,3) {\footnotesize $g(x^1_j)$};
\node[draw,circle,fill=white] at (1,0) {\footnotesize $g(x^4_k)$};
\node[draw,circle,fill=white] at (1,1) {\footnotesize $g(x^3_k)$};
\node[draw,circle,fill=white] at (1,2) {\footnotesize $g(x^2_k)$};
\node[draw,circle,fill=white] at (1,3) {\footnotesize $g(x^1_k)$};

\node[above] at (0,3.5) {$\mathbf{H_j}$};
\node[above] at (1,3.5) {$\mathbf{H_k}$};
\node[right] at (1.5,0) {$\mathbf{\ell^4}$};
\node[right] at (1.5,1) {$\mathbf{\ell^3}$};
\node[right] at (1.5,2) {$\mathbf{\ell^2}$};
\node[right] at (1.5,3) {$\mathbf{\ell^1}$};
\node[right] at (0,2.5) {$\mathbf{\neq}$};
\node[left] at (1,0.5) {$\mathbf{\neq}$};
\node[above] at (0.5,3) {$\mathbf{\neq}$};
\node[above] at (0.5,2) {$\mathbf{\neq}$};
\node[above] at (0.5,1) {$\mathbf{\neq}$};
\node[above] at (0.5,0) {$\mathbf{\neq}$};

\end{tikzpicture}
\caption{A forbidden configuration in WTT game forms as in Lemma \ref{forbidden}.\label{fig1}}
\end{figure}
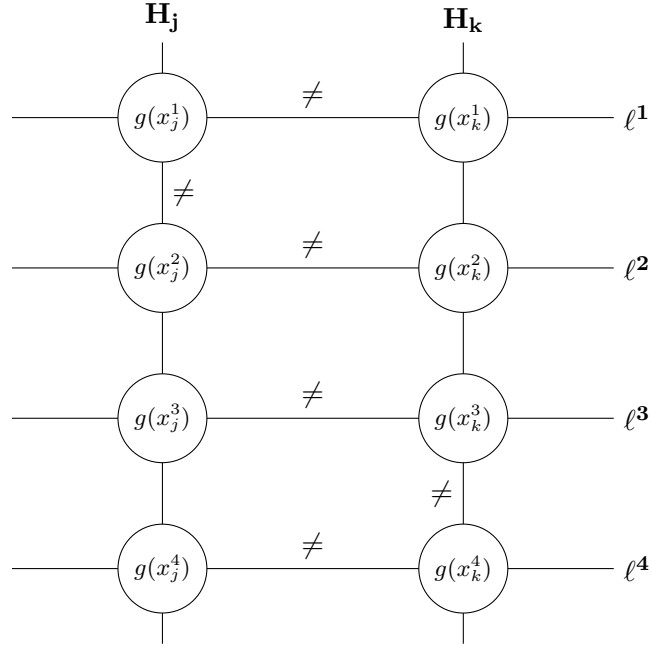

\proof
Using the inequalities (1) and (3) for the quadruple $(x_j^1, x_j^2, x_k^1, x_k^2)$ we either get a direct contradiction to the WTT property of $g$ or we get $g(x_k^1) = g(x_k^2) = a_k$ for some outcome $a_k \in A$. Similarly using (2) and (3) for the quadruple $(x_j^3, x_j^4, x_k^3, x_k^4)$ we either get a direct contradiction to the WTT property of $g$ or we get $g(x_j^3) = g(x_j^4) = a_j$ for some outcome $a_j \in A$. So let us assume the latter in both cases. Now using (1) we get that one of $g(x_j^1), g(x_j^2)$ must differ from $a_j$, so let us denote the differing strategy profile by $x_j^u$, for $u \in \{1,2\}$. Similarly, using (2) we get $g(x_k^3) \neq a_k$ or $g(x_k^4) \neq a_k$ and let us denote the differing strategy profile by $x_k^v$, for $v \in \{3,4\}$. This altogether implies that the quadruple $(x_j^u, x_j^v, x_k^u, x_k^v)$ contradicts the WTT property of $g$, see Figure \ref{fig1}.

Note that if some of the lines $\ell^i, 1\leq i\leq 4$ coincide (of course by the assumptions $\ell_1$ must differ from $\ell_2$ and $\ell_3$ must differ from $\ell_4$), the proof becomes even simpler. If $\ell_1 = \ell_3$ and $\ell_2 = \ell_4$ (or $\ell_1 = \ell_4$ and $\ell_2 = \ell_3$) then the four intersections immediately give a contradiction to the WTT property. If $\ell_1 = \ell_3$ and $\ell_2 \neq \ell_4$ then the proof above goes through for $u=2$ and $v=4$. Symmetrically, if $\ell_1 \neq \ell_3$ and $\ell_2 = \ell_4$ then the proof above goes through for $u=1$ and $v=3$.
\qed

Let us now define a notation for special hyperplane partitions and
state and prove a key property of these partitions.

\begin{definition}
\label{d-hyperplane-partition}
Let $g$ be an $n$-person game form and $i$ an arbitrary direction (player). Let $H_j$ and $H_k$ be two distinct parallel hyperplanes perpendicular to direction $i$. For an arbitrary line $\ell$ in direction $i$ let us denote by $x_j^{\ell}$ and $x_k^{\ell}$ the strategy profiles at the intersections of line $\ell$ with hyperplanes $H_j$ and $H_k$ respectively. We define a partition of $H_j$ into $H_j^=(k)$ and $H_j^{\neq}(k)$ as follows:
\[ H_j^=(k) = \{ x_j^{\ell} \ST g(x_j^{\ell}) = g(x_k^{\ell})\} \;\;\; \mbox{ and } \;\;\;
H_j^{\neq}(k) = \{ x_j^{\ell} \ST g(x_j^{\ell}) \neq g(x_k^{\ell})\} \]
\end{definition}

\begin{lemma}
\label{constant-region}
Let $g$ be an $n$-person WTT game form, $i$ an arbitrary direction (player), and $H_j, H_k$ be two distinct parallel hyperplanes perpendicular to direction $i$. Then $H_j^{\neq}(k)$ is a constant region or $H_k^{\neq}(j)$ is a constant region (or both are constant regions).
\end{lemma}

\proof
Assume by contradiction that neither $H_j^{\neq}(k)$ nor $H_k^{\neq}(j)$ is a constant region. This means that there exist four lines $\ell^i, 1\leq i\leq 4$ (not necessarily all distinct) in direction $i$ intersecting $H_j^{\neq}(k)$ in strategy profiles $x_j^i, 1\leq i\leq 4$ and intersecting $H_k^{\neq}(j)$ in strategy profiles $x_k^i, 1\leq i\leq 4$  such that
\begin{description}
	\item[(1)] $g(x_j^1) \neq g(x_j^2)$ (because $H_j^{\neq}(k)$ contains two distinct outcomes),
	\item[(2)] $g(x_k^3) \neq g(x_k^4)$ (because $H_k^{\neq}(j)$ contains two distinct outcomes),
	\item[(3)] $g(x_j^i) \neq g(x_k^i), 1\leq i\leq 4$ (by the definition of $H_j^{\neq}(k)$ and $H_k^{\neq}(j)$).
\end{description}
The four lines $\ell^i, 1\leq i\leq 4$ and their intersections with hyperplanes $H_j, H_k$ obviously fulfil the assumptions of Lemma~\ref{forbidden} and hence $g$ is not WTT which is a contradiction.
\qed

\begin{definition}
\label{hyperplane-domination}
Let $g$ be an $n$-person WTT game form, $i$ an arbitrary direction (player), and $H_j, H_k$ be two distinct parallel hyperplanes perpendicular to direction $i$. If $H_j^{\neq}(k)$ is a constant $c$ region for some outcome $c$, then we say that $H_j$ dominates $H_k$ by $c$ and denote this fact by $H_j \stackrel{c}{\longrightarrow} H_k$. If $H_j \stackrel{c}{\longrightarrow} H_k$ and there exists no outcome $d$ such that $H_k \stackrel{d}{\longrightarrow} H_j$ then we say that $H_j$ strictly dominates $H_k$ by $c$ and write $H_j \stackrel{c}{\Longrightarrow} H_k$.
\end{definition}

Note that $H_j \stackrel{c}{\Longrightarrow} H_k$ implies that $H_k^{\neq}(j)$ is not a constant region. Using the just defined notation, Lemma~\ref{constant-region} and the fact that we have no two identical parallel hyperplanes by Remark \ref{basic-assumptions} implies the following easy corollary.

\begin{corollary}
\label{graph-options}
Let $g$ be an $n$-person WTT game form, $i$ an arbitrary direction (player), and $H_j, H_k$ be two distinct parallel hyperplanes perpendicular to direction $i$. Then exactly one of the following three conditions is true
\begin{enumerate}
\item $H_j \stackrel{c}{\Longrightarrow} H_k$ for some outcome $c$ and there exist two distinct outcomes $a \neq b$ in $H_k^{\neq}(j)$ (which are both different from $c$)
\item $H_k \stackrel{d}{\Longrightarrow} H_j$ for some outcome $d$ and there exist two distinct outcomes $a \neq b$ in $H_j^{\neq}(k)$ (which are both different from $d$)
\item $H_j \stackrel{c}{\longrightarrow} H_k$ and $H_k \stackrel{d}{\longrightarrow} H_j$ for some outcomes $c \neq d$.
\end{enumerate}
\end{corollary}

\begin{remark}
\label{WTT-characterization}
It should be noted that Corollary~\ref{graph-options} gives a complete characterization of WTT game forms. Namely, a game form is WTT if and only if any pair of parallel hyperplanes fulfills exactly one of the three properties specified in Corollary~\ref{graph-options}. The left to right implications is proved in Lemma~\ref{constant-region} while the reverse implication is trivial.
\end{remark}

Now we shall show that a hyperplane cannot strictly dominate two other hyperplanes by two different outcomes.

\begin{lemma}
\label{strict-domination}
Let $g$ be an $n$-person WTT game form, $i$ an arbitrary direction (player), and $H_{\ell}, H_j, H_k$ three distinct parallel hyperplanes perpendicular to direction $i$ such that $H_{\ell} \stackrel{c}{\Longrightarrow} H_j$ and $H_{\ell} \stackrel{d}{\Longrightarrow} H_k$ for some outcomes $c$ and $d$. Then $c=d$.
\end{lemma}

\begin{figure}[ht]
\centering
\begin{tikzpicture}[x={(4cm,0cm)},y={(0cm,2cm)}]

\draw[thin,step=1] (-0.5,-0.5) grid (2.5,3.5);

\node[draw,circle,fill=white] at (0,0) {\footnotesize $\color{blue}d$};
\node[draw,circle,fill=white] at (0,1) {\footnotesize $\color{blue}d$};
\node[draw,circle,fill=white] at (0,2) {\footnotesize $a_2$};
\node[draw,circle,fill=white] at (0,3) {\footnotesize $a_1$};
\node[draw,circle,fill=white] at (1,0) {\footnotesize $d$};
\node[draw,circle,fill=white] at (1,1) {\footnotesize $d$};
\node[draw,circle,fill=white] at (1,2) {\footnotesize $c$};
\node[draw,circle,fill=white] at (1,3) {\footnotesize $c$};
\node[draw,circle,fill=white] at (2,0) {\footnotesize $b_2$};
\node[draw,circle,fill=white] at (2,1) {\footnotesize $b_1$};
\node[draw,circle,fill=white] at (2,2) {\footnotesize $\color{red}c$};
\node[draw,circle,fill=white] at (2,3) {\footnotesize $\color{red}c$};

\node[above] at (0,3.5) {$\mathbf{H_j}$};
\node[above] at (1,3.5) {$\mathbf{H_\ell}$};
\node[above] at (2,3.5) {$\mathbf{H_k}$};

\node[above] at (1.5,3.5) {$\color{red}\stackrel{d}{\Longrightarrow}$};
\node[above] at (0.5,3.5) {$\color{blue}\stackrel{c}{\Longleftarrow}$};

\node[right] at (2.5,0) {$\mathbf{\ell^4}$};
\node[right] at (2.5,1) {$\mathbf{\ell^3}$};
\node[right] at (2.5,2) {$\mathbf{\ell^2}$};
\node[right] at (2.5,3) {$\mathbf{\ell^1}$};

\node[right] at (2.6,2.5) {$\left.\rule{0mm}{2cm}\right\} H^{\neq}_\ell(j) \subseteq {\color{red}H^{=}_\ell(k)}$};
\node[right] at (2.6,0.5) {$\left.\rule{0mm}{2cm}\right\} H^{\neq}_\ell(k) \subseteq {\color{blue}H^{=}_\ell(j)}$};

\node[right] at (0,2.5) {$\mathbf{\neq}$};
\node[left] at (2,0.5) {$\mathbf{\neq}$};
\node[above] at (0.5,3) {$\mathbf{\neq}$};
\node[above] at (0.5,2) {$\mathbf{\neq}$};
\node[above] at (1.5,1) {$\mathbf{\neq}$};
\node[above] at (1.5,0) {$\mathbf{\neq}$};

\node[above] at (1.5,3) {$\color{red}\mathbf{=}$};
\node[above] at (1.5,2) {$\color{red}\mathbf{=}$};

\node[above] at (0.5,1) {$\color{blue}\mathbf{=}$};
\node[above] at (0.5,0) {$\color{blue}\mathbf{=}$};

\end{tikzpicture}
\caption{Strict dominance by two different outcomes $c\neq d$ leads to the forbidden configuration as in Lemma \ref{forbidden}.\label{fig2}}
\end{figure}
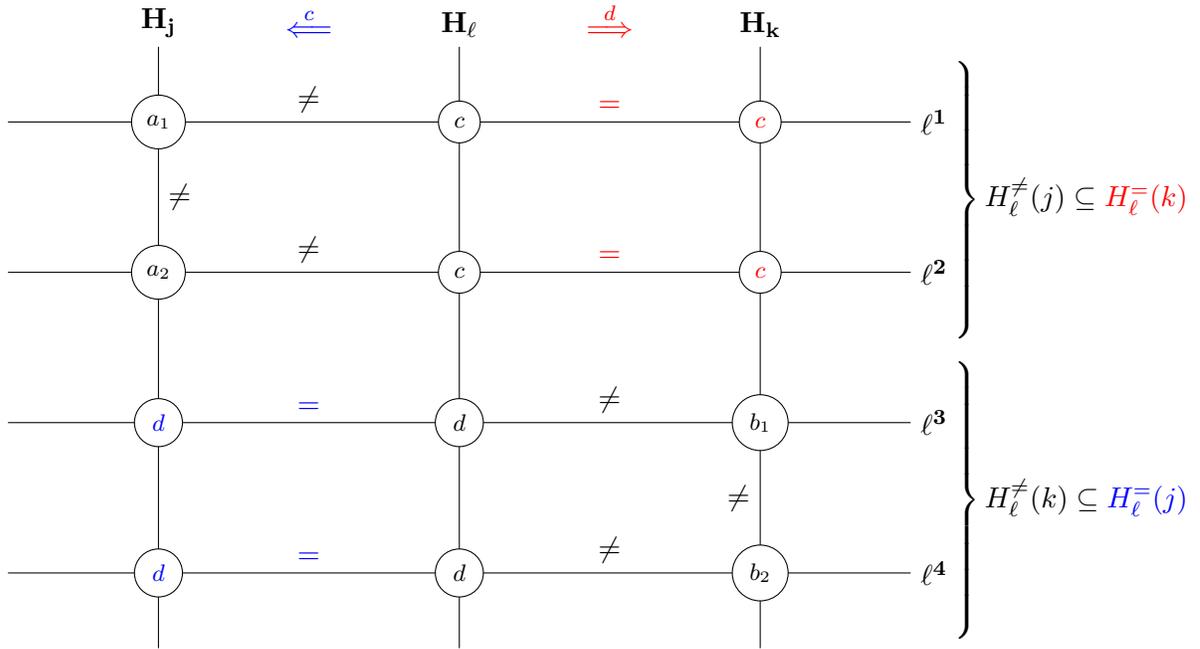

\proof
Assume by contradiction that $c \neq d$. Then, hyperplane $H_{\ell}$ can be partitioned into three disjoint regions, namely the constant $c$ region $H_{\ell}^{\neq}(j)$, a constant $d$ region $H_{\ell}^{\neq}(k)$, and region $H_{\ell}^=(j)\cap H_{\ell}^=(k)$ where the outcomes are the same in all three hyperplanes for any perpendicular line in direction $i$. This in particular implies that $H_{\ell}^{\neq}(j) \subseteq H_{\ell}^=(k)$ and $H_{\ell}^{\neq}(k) \subseteq H_{\ell}^=(j)$, see Figure \ref{fig2}.

Now $H_\ell \stackrel{c}{\Longrightarrow} H_j$ implies that there exist two distinct lines $\ell_1, \ell_2$ in direction $i$ intersecting $H_{\ell}^{\neq}(j)$ (and thus also $H_{\ell}^=(k)$) in profiles $x_{\ell}^1, x_{\ell}^2$ for which $g(x_{\ell}^1) = g(x_{\ell}^2) = c$, intersecting $H_j^{\neq}(\ell)$ in profiles $x_j^1, x_j^2$ for which $g(x_j^1) = a_1 \neq a_2 = g(x_j^2)$ (here we use the strict domination), and intersecting $H_k^=(\ell)$ in profiles $x_k^1, x_k^2$ for which $g(x_k^1) = g(x_k^2) = c$. Note that the outcomes $a_1, a_2, c$ are pairwise distinct.

Similarly, $H_{\ell} \stackrel{d}{\Longrightarrow} H_k$ implies that there exist two distinct lines $\ell_3, \ell_4$ in direction $i$ intersecting $H_{\ell}^{\neq}(k)$ (and thus also $H_{\ell}^=(j)$) in profiles $x_{\ell}^3, x_{\ell}^4$ for which $g(x_{\ell}^3) = g(x_{\ell}^4) = d$, intersecting $H_j^=(\ell)$ in profiles $x_j^3, x_j^4$ for which $g(x_j^3) = g(x_j^4) = d$, and intersecting $H_k^{\neq}(\ell)$ in profiles $x_k^3, x_k^4$ for which $g(x_k^3) = b_1$, $g(x_k^3) = b_2$ for $b_1, b_2, d$ pairwise distinct.

It is easy to check that the four pairwise distinct lines $\ell^i, 1\leq i\leq 4$ and their intersections with hyperplanes $H_j, H_k$ fulfil the assumptions of Lemma~\ref{forbidden} and hence $g$ is not WTT which is a contradiction.
\qed

Lemma~\ref{strict-domination} allows us to associate a unique outcome to every hyperplane that strictly dominates at least one other hyperplane.

\begin{definition}
\label{proper-outcome}
Let $g$ be an $n$-person WTT game form, $i$ an arbitrary direction (player), and $H_j$ be a hyperplane perpendicular to direction $i$. If there exists a hyperplane $H_k$ parallel to $H_j$ and an output $c$ such that $H_j \stackrel{c}{\Longrightarrow} H_k$ then we call $c$ the {\em proper outcome} of $H_j$.
\end{definition}

Note that a hyperplane $H_j$ that does not strictly dominate any other hyperplane must have the property (by Corollary~\ref{graph-options}), that it is dominated (non-strictly or strictly) by every other hyperplane parallel to $H_j$. We shall call such a hyperplane for which no proper outcome is defined a sink hyperplane.

\begin{definition}
\label{sink-hyperplane}
Let $g$ be an $n$-person WTT game form, $i$ an arbitrary direction (player), and $H_j$ be a hyperplane perpendicular to direction $i$. We shall call $H_j$ a {\em sink hyperplane} if for every hyperplane $H_k$, $k \neq j$, perpendicular to direction $i$ there exists an outcome $c_k$ such that $H_k \stackrel{c_k}{\longrightarrow} H_j$. If there exist no sink hyperplane in any direction then $g$ is called a {\em no-sink game form}.
\end{definition}

Definitions~\ref{proper-outcome} and~\ref{sink-hyperplane} allow us to formulate the following simple corollary.

\begin{corollary}
\label{sink-or-proper}
Let $g$ be an $n$-person WTT game form, $i$ an arbitrary direction (player), and $H_j$ be a hyperplane perpendicular to direction $i$. Then either $H_j$ has a proper outcome or it is a sink hyperplane.
\end{corollary}

\begin{figure}[ht]
\centering
\begin{tikzpicture}[z={(0.45cm,0.6cm)},x={(1.2cm,0cm)},y={(0cm,1cm)},scale=2]

\foreach \x in {1,2,3} \foreach \y in {1,2,3}
	\draw[gray,thin] (\x,\y,1) -- (\x,\y,3);

\foreach \x in {1,2,3} \foreach \z in {1,2,3}
	\draw[gray,thin] (\x,1,\z) -- (\x,3,\z);
	
\foreach \y in {1,2,3} \foreach \z in {1,2,3}
	\draw[gray,thin] (1,\y,\z) -- (3,\y,\z);

\node[circle,fill=white] at (1,1,1) {$a$};
\node[circle,fill=white] at (1,2,1) {$a$};
\node[circle,fill=white] at (1,3,1) {$a$};
\node[circle,fill=white] at (1,1,2) {$a$};
\node[circle,fill=white] at (1,2,2) {$a$};
\node[circle,fill=white] at (1,3,2) {$a$};
\node[circle,fill=white] at (1,1,3) {$c$};
\node[circle,fill=white] at (1,2,3) {$a$};
\node[circle,fill=white] at (1,3,3) {$a$};

\node[circle,fill=white] at (2,1,1) {$a$};
\node[circle,fill=white] at (2,2,1) {$a$};
\node[circle,fill=white] at (2,3,1) {$a$};
\node[circle,fill=white] at (2,1,2) {$b$};
\node[circle,fill=white] at (2,2,2) {$b$};
\node[circle,fill=white] at (2,3,2) {$a$};
\node[circle,fill=white] at (2,1,3) {$b$};
\node[circle,fill=white] at (2,2,3) {$b$};
\node[circle,fill=white] at (2,3,3) {$a$};

\node[circle,fill=white] at (3,1,1) {$c$};
\node[circle,fill=white] at (3,2,1) {$a$};
\node[circle,fill=white] at (3,3,1) {$a$};
\node[circle,fill=white] at (3,1,2) {$b$};
\node[circle,fill=white] at (3,2,2) {$b$};
\node[circle,fill=white] at (3,3,2) {$a$};
\node[circle,fill=white] at (3,1,3) {$c$};
\node[circle,fill=white] at (3,2,3) {$b$};
\node[circle,fill=white] at (3,3,3) {$c$};

\shade[right color=blue,left color=white,top color=blue,bottom color=white,opacity=0.2] (1,3,3) -- (1,3,4) -- (1,2.5,4) -- cycle;
\node[above right] at (1,3,4.25) {$\color{blue}\mathbf{H_1}$};
\shade[right color=blue,left color=blue,opacity=0.2] (1,1,1) -- (1,3,1) -- (1,3,3) -- (1,1,3) -- cycle;

\shade[right color=red,left color=white,top color=red,bottom color=white,opacity=0.2] (2,3,3) -- (2,3,4) -- (2,2.5,4) -- cycle;
\node[above right] at (2,3,4) {$\color{red}\mathbf{H_2}$};
\shade[right color=red,left color=red,opacity=0.2] (2,1,1) -- (2,3,1) -- (2,3,3) -- (2,1,3) -- cycle;

\shade[right color=green,left color=white,top color=green,bottom color=white,opacity=0.2] (3,3,3) -- (3,3,4) -- (3,2.5,4) -- cycle;
\node[above right] at (3,3,4.25) {$\color{green}\mathbf{H_3}$};
\shade[right color=green,left color=green,opacity=0.2] (3,1,1) -- (3,3,1) -- (3,3,3) -- (3,1,3) -- cycle;

\node[above right] at (1.5,3,4) {$\stackrel{b}{\Longleftarrow}$};
\node[above right] at (2.5,3,4) {$\stackrel{c}{\Longleftarrow}$};
\node[above] at (2,3,4.5) {$\stackrel{a}{\Longrightarrow}$};

\end{tikzpicture}
\caption{First no-sink 3D example.\label{f-nosink-3D-1}}
\end{figure}
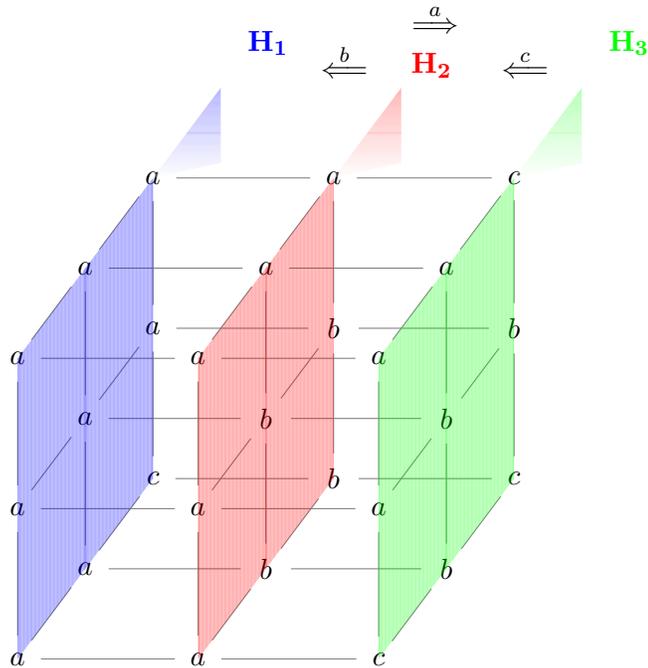

We will now study no-sink WTT game forms. In such a game form every hyperplane strictly dominates at least one of its parallel hyperplanes, which in turn implies that there must be a cycle (or several cycles) of strong dominance relations among all hyperplanes in every direction. Note that such game forms exist, consider for instance the following 2-person game form
\begin{equation}
\label{2D-no-sink}
\left[ {\begin{array}{ccc}
a & a & c \\
a & b & b \\
c & b & c
\end{array} } \right]
\end{equation}
in which $R_1 \stackrel{a}{\Longrightarrow} R_3$, $R_3 \stackrel{c}{\Longrightarrow} R_2$, and $R_2 \stackrel{b}{\Longrightarrow} R_1$, where $R_i$ is the $i$-th row of the game form. The same cycle of strict dominance holds among the columns $C_1$, $C_2$, and $C_3$ due to the symmetry w.r.t. the main diagonal. Generating an example for three players is much more difficult to do by hand, but two such game forms were computer generated using a code of V.Oudalov. Each of them is a $3 \times 3 \times 3$ array displayed below as a set of three 2-dimensional $3 \times 3$ subarrays (hyperplanes). The first example is shown in Figure \ref{f-nosink-3D-1}.

Clearly $H_1 \stackrel{a}{\Longrightarrow} H_3$, $H_3 \stackrel{c}{\Longrightarrow} H_2$, and $H_2 \stackrel{b}{\Longrightarrow} H_1$. Note that $H_3$ is exactly the two dimensional game form from~(\ref{2D-no-sink}) and the relations from there carry over to the row hyperplanes $R_i$, $1 \leq i \leq 3$ (in particular $R_1 \stackrel{a}{\Longrightarrow} R_3$, $R_3 \stackrel{c}{\Longrightarrow} R_2$, and $R_2 \stackrel{b}{\Longrightarrow} R_1$) and due to symmetry also to the column hyperplanes $C_i$, $1 \leq i \leq 3$.

\begin{figure}[ht]
\centering
\begin{tikzpicture}[z={(0.45cm,0.6cm)},x={(1.2cm,0cm)},y={(0cm,1cm)},scale=2]

\foreach \x in {1,2,3} \foreach \y in {1,2,3}
	\draw[gray,thin] (\x,\y,1) -- (\x,\y,3);

\foreach \x in {1,2,3} \foreach \z in {1,2,3}
	\draw[gray,thin] (\x,1,\z) -- (\x,3,\z);
	
\foreach \y in {1,2,3} \foreach \z in {1,2,3}
	\draw[gray,thin] (1,\y,\z) -- (3,\y,\z);

\node[circle,fill=white] at (1,1,1) {$a$};
\node[circle,fill=white] at (1,2,1) {$a$};
\node[circle,fill=white] at (1,3,1) {$a$};
\node[circle,fill=white] at (1,1,2) {$a$};
\node[circle,fill=white] at (1,2,2) {$a$};
\node[circle,fill=white] at (1,3,2) {$a$};
\node[circle,fill=white] at (1,1,3) {$c$};
\node[circle,fill=white] at (1,2,3) {$a$};
\node[circle,fill=white] at (1,3,3) {$a$};

\node[circle,fill=white] at (2,1,1) {$a$};
\node[circle,fill=white] at (2,2,1) {$a$};
\node[circle,fill=white] at (2,3,1) {$a$};
\node[circle,fill=white] at (2,1,2) {$b$};
\node[circle,fill=white] at (2,2,2) {$b$};
\node[circle,fill=white] at (2,3,2) {$a$};
\node[circle,fill=white] at (2,1,3) {$b$};
\node[circle,fill=white] at (2,2,3) {$b$};
\node[circle,fill=white] at (2,3,3) {$a$};

\node[circle,fill=white] at (3,1,1) {$c$};
\node[circle,fill=white] at (3,2,1) {$c$};
\node[circle,fill=white] at (3,3,1) {$c$};
\node[circle,fill=white] at (3,1,2) {$c$};
\node[circle,fill=white] at (3,2,2) {$b$};
\node[circle,fill=white] at (3,3,2) {$c$};
\node[circle,fill=white] at (3,1,3) {$c$};
\node[circle,fill=white] at (3,2,3) {$c$};
\node[circle,fill=white] at (3,3,3) {$c$};

\shade[right color=blue,left color=white,top color=blue,bottom color=white,opacity=0.2] (1,3,3) -- (1,3,4) -- (1,2.5,4) -- cycle;
\node[above right] at (1,3,4.25) {$\color{blue}\mathbf{H_1}$};
\shade[right color=blue,left color=blue,opacity=0.2] (1,1,1) -- (1,3,1) -- (1,3,3) -- (1,1,3) -- cycle;

\shade[right color=red,left color=white,top color=red,bottom color=white,opacity=0.2] (2,3,3) -- (2,3,4) -- (2,2.5,4) -- cycle;
\node[above right] at (2,3,4) {$\color{red}\mathbf{H_2}$};
\shade[right color=red,left color=red,opacity=0.2] (2,1,1) -- (2,3,1) -- (2,3,3) -- (2,1,3) -- cycle;

\shade[right color=green,left color=white,top color=green,bottom color=white,opacity=0.2] (3,3,3) -- (3,3,4) -- (3,2.5,4) -- cycle;
\node[above right] at (3,3,4.25) {$\color{green}\mathbf{H_3}$};
\shade[right color=green,left color=green,opacity=0.2] (3,1,1) -- (3,3,1) -- (3,3,3) -- (3,1,3) -- cycle;

\node[above right] at (1.5,3,4) {$\stackrel{b}{\Longleftarrow}$};
\node[above right] at (2.5,3,4) {$\stackrel{c}{\Longleftarrow}$};
\node[above] at (2,3,4.5) {$\stackrel{a}{\Longrightarrow}$};

\end{tikzpicture}
\caption{Second no-sink 3D example.\label{f-no-sink-3D-2}}
\end{figure}
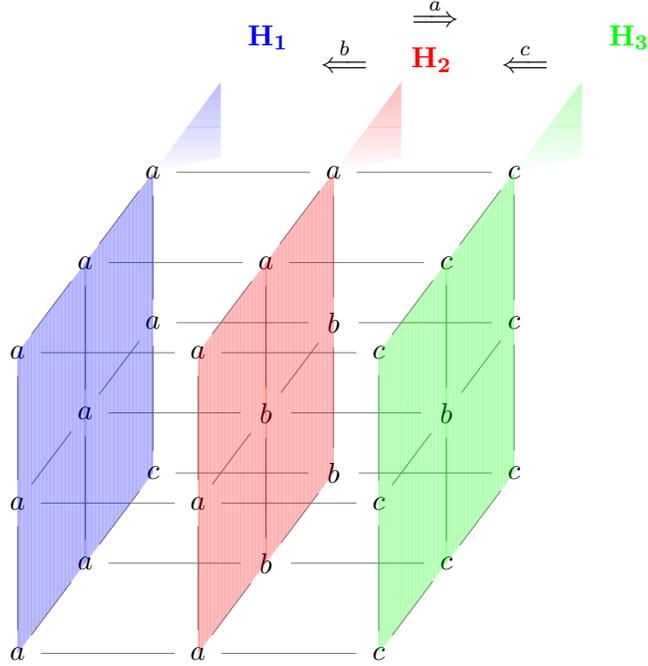

The second example we have is shown in Figure \ref{f-no-sink-3D-2}.
Note that only $H_3$ is different. We leave the detection of the three strong dominance cycles in all three directions to the reader. We conjecture, that no-sink game forms exist for any $n$, not just for $n=2$ and $n=3$ as displayed above. Before we start to study the properties of no-sink WTT game forms let us introduce two definitions.

\begin{definition}
\label{point-direction-hyperplane}
Let $g$ be an $n$-person WTT game form, $x = (x_1, \ldots ,x_n)$ a strategy profile, and $i$ an arbitrary direction (player). Then $H_i^x$ denotes the hyperplane perpendicular to direction $i$ containing the profile $x$. That is, $H_i^x$ consists of all strategy profiles in $g$ which have the $i$-th coordinate equal to $x_i$.
\end{definition}

\begin{definition}
\label{k-box}
Let $g$ be an $n$-person no-sink WTT game form. We say that $g$ contains a {\em $k$-box} if there exist two strategy profiles $x = (x_1, \ldots ,x_n)$ and $y = (y_1, \ldots ,y_n)$ such that:
\begin{enumerate}
\item $g(x) \neq g(y)$,
\item $x$ and $y$ differ in exactly $k$ coordinates $i_1, \ldots ,i_k$ (so $x$ and $y$ span a k-dimensional subgame form of $2^k$ strategy profiles), and
\item for every $1 \leq j \leq k$ it holds that $g(x)$ is not the proper outcome of $H_{i_j}^x$ and $g(y)$ is not the proper outcome of $H_{i_j}^y$.
\end{enumerate}
\end{definition}

Now we are ready to state several properties of no-sink WTT game forms.

\begin{lemma}
\label{box-induction}
Let $g$ be an $n$-person no-sink WTT game form which contains a $k$-box for some $1 \leq k \leq n$. Then $g$ contains a $(k-1)$-box or a $1$-box.
\end{lemma}

\proof
Let us assume without a loss of generality that profiles $x$ and $y$ spanning the $k$-box differ in the first $k$ coordinates (if not we can permute the coordinates), i.e. $x = (x_1, \ldots ,x_n)$ and $y = (y_1, \ldots ,y_n)$ where $x_i \neq y_i$ for $1 \leq i \leq k$ and $x_i = y_i$ for $k+1 \leq i \leq n$. Now consider direction $1$ and strategy profiles $x' = (y_1,x_2 \ldots ,x_n)$ and $y' = (x_1, y_2 \ldots ,y_n)$. Notice that the pairs $x,x'$ and $y,y'$ lie on lines in direction $1$ while the pairs $x,y'$ and $x',y$ belong to two hyperplanes perpendicular to direction $1$ (namely $H_1^x$ and $H_1^y$). Thus the rectangle $x,x',y,y'$ must contain a constant line due to the WTT property. Now we have four possibilities:
\begin{enumerate}
\item if $g(x') = g(y)$ then $x,x'$ span a $1$-box,
\item if $g(x') = g(x)$ then $y,x'$ span a $(k-1)$-box,
\item if $g(y') = g(y)$ then $x,y'$ span a $(k-1)$-box, and
\item if $g(y') = g(x)$ then $y,y'$ span a $1$-box.
\end{enumerate}
In all four cases the three properties defining the specified $1$-box or $(k-1)$-box follow easily from the properties of the $k$-box. In particular, in the first case $g(x) \neq g(x') = g(y)$, $x$ and $x'$ differ only in the first coordinate, and a pair of parallel hyperplanes $H_1^{x}$ and $H_1^{x'} = H_1^{y}$  fulfills the third property. In the second case $g(y) \neq g(x') = g(x)$, $y$ and $x'$ differ in exactly $(k-1)$ coordinates (namely $2, \ldots ,k$), and $(k-1)$ pairs of parallel hyperplanes $H_i^{y}$ and $H_i^{x'} = H_i^{x}$ for $2 \leq i \leq k$ fulfill the third property. The third and fourth case are symmetric.
\qed

\begin{lemma}
\label{no-1-box}
Let $g$ be an $n$-person no-sink WTT game form. Then $g$ contains no $1$-box.
\end{lemma}

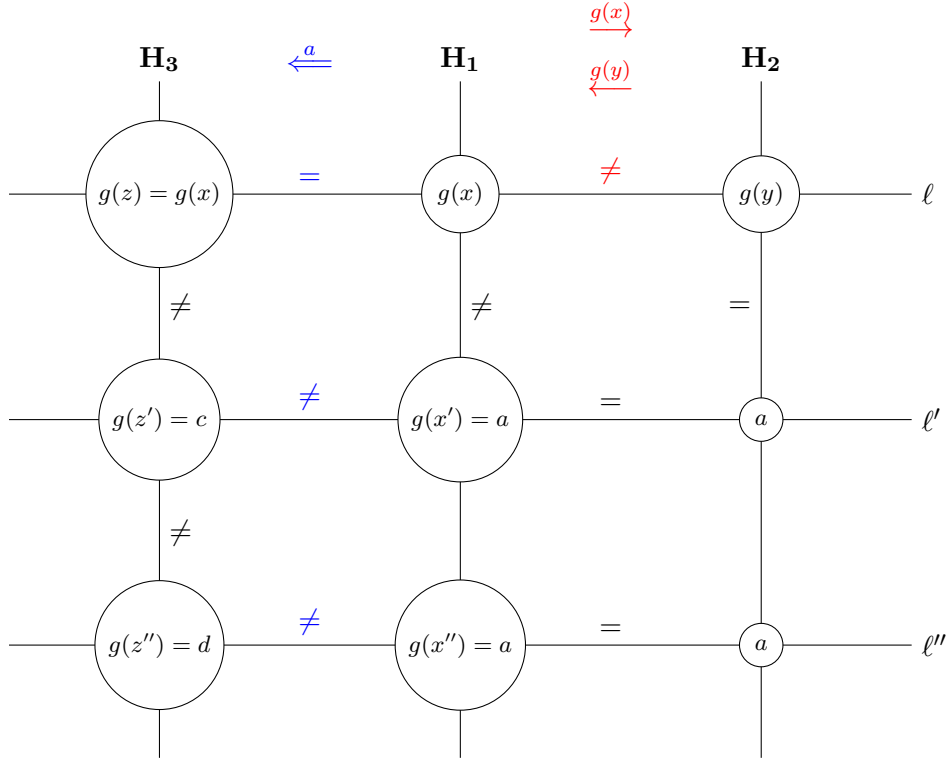
\begin{figure}[ht]
\centering
\begin{tikzpicture}[x={(4cm,0cm)},y={(0cm,3cm)}]

\draw[thin,step=1] (-0.5,-0.5) grid (2.5,2.5);

\node[draw,circle,fill=white] at (0,0) {\footnotesize $g(z'')=d$};
\node[draw,circle,fill=white] at (0,1) {\footnotesize $g(z')=c$};
\node[draw,circle,fill=white] at (0,2) {\footnotesize $g(z)=g(x)$};
\node[draw,circle,fill=white] at (1,0) {\footnotesize $g(x'')=a$};
\node[draw,circle,fill=white] at (1,1) {\footnotesize $g(x')=a$};
\node[draw,circle,fill=white] at (1,2) {\footnotesize $g(x)$};
\node[draw,circle,fill=white] at (2,0) {\footnotesize $a$};
\node[draw,circle,fill=white] at (2,1) {\footnotesize $a$};
\node[draw,circle,fill=white] at (2,2) {\footnotesize $g(y)$};

\node[above] at (0,2.5) {$\mathbf{H_3}$};
\node[above] at (1,2.5) {$\mathbf{H_1}$};
\node[above] at (2,2.5) {$\mathbf{H_2}$};

\node[above] at (1.5,2.65) {$\color{red}\stackrel{g(x)}{\longrightarrow}$};
\node[below] at (1.5,2.65) {$\color{red}\stackrel{g(y)}{\longleftarrow}$};
\node[above] at (0.5,2.5) {$\color{blue}\stackrel{a}{\Longleftarrow}$};

\node[right] at (2.5,0) {$\mathbf{\ell''}$};
\node[right] at (2.5,1) {$\mathbf{\ell'}$};
\node[right] at (2.5,2) {$\mathbf{\ell}$};

\node[above] at (1.5,0) {$\mathbf{=}$};
\node[above] at (1.5,1) {$\mathbf{=}$};
\node[above] at (1.5,2) {$\color{red}\mathbf{\neq}$};

\node[above] at (0.5,0) {$\color{blue}\mathbf{\neq}$};
\node[above] at (0.5,1) {$\color{blue}\mathbf{\neq}$};
\node[above] at (0.5,2) {$\color{blue}\mathbf{=}$};

\node[right] at (0,0.5) {$\mathbf{\neq}$};
\node[right] at (0,1.5) {$\mathbf{\neq}$};

\node[right] at (1,1.5) {$\mathbf{\neq}$};

\node[left] at (2,1.5) {$\mathbf{=}$};

\end{tikzpicture}
\caption{Configuration showing that no $1$-box can exists in a WTT game form.\label{fig-no-1-box}}
\end{figure}

\proof
By contradiction let $x$ and $y$ be two profiles spanning a $1$-box, i.e. such that
\begin{enumerate}
\item $g(x) \neq g(y)$,
\item $x$ and $y$ differ in exactly $1$ coordinate, i.e. lie on a line $\ell$ in some direction $i$, and
\item $g(x)$ is not the proper outcome of $H_{i}^x$ and $g(y)$ is not the proper outcome of $H_{i}^y$.
\end{enumerate}
Let us denote the proper outcome of $H_1 = H_{i}^x$ by $a \neq g(x)$ and the proper outcome of $H_2 = H_{i}^y$ by $b \neq g(y)$. Since $a \neq g(x) \neq g(y) \neq b$ we can have neither $H_1 \stackrel{a}{\Longrightarrow} H_2$ nor $H_2 \stackrel{b}{\Longrightarrow} H_1$. Consequently, $H_1$ cannot strictly dominate $H_2$ by Lemma \ref{strict-domination}, since we assumed $a$ to be its proper outcome.  Similarly $H_2$ cannot strictly dominate $H_1$. Therefore, by Corollary \ref{graph-options} we must have both $H_1\stackrel{g(x)}{\longrightarrow} H_2$ and $H_2\stackrel{g(y)}{\longrightarrow} H_1$, see Figure \ref{fig-no-1-box}.

Hence, there must exist a third hyperplane $H_3$ perpendicular to direction $i$ and distinct from $H_1, H_2$ such that $H_1 \stackrel{a}{\Longrightarrow} H_3$. This implies that line $\ell$ intersects $H_3$ in some strategy profile $z$ with $g(z) = g(x)$, because $x$ belongs to $H_1^=(3)$. Moreover, there must exist two distinct lines $\ell'$ and $\ell''$ in direction $i$ intersecting $H_1^{\neq}(3)$ in profiles $x'$ and $x''$ with $g(x') = g(x'') = a$ and intersecting $H_3^{\neq}(1)$ in profiles $z'$ and $z''$ where $g(z') = c$, $g(z'') = d$ with $a,c,d$ pairwise distinct. One of $c,d$ must be different from $g(x)$ so let us assume $c \neq g(x) = g(z)$. Let $\ell'$ and $\ell''$ intersect $H_2$ in profiles $y'$ and $y''$. Because $H_1 \stackrel{g(x)}{\longrightarrow} H_2$, lines $\ell'$ and $\ell''$ intersect $H_1^=(2)$ and so $g(y') = g(y'') = a$. Now consider the rectangle $y,z,z',y'$ on lines $\ell, \ell'$ and in hyperplanes $H_2,H_3$. We have $g(y) \neq g(z) = g(x)$, $g(z) \neq g(z') = c$, and $c = g(z') \neq g(y') = a$. Therefore, WTT property implies $g(y) = g(y') = a$. However, now the quadruple $y',y'',z',z''$ with $g(y') = g(y'') = g(y) = a$, $g(z') = c$, $g(z'') = d$ with $a,c,d$ pairwise distinct implies $H_2 \stackrel{g(y)}{\Longrightarrow} H_3$ which contradicts the fact that the proper outcome of $H_2$ is $b \neq g(y)$.
\qed

Lemma~\ref{box-induction} and Lemma~\ref{no-1-box} of course have an obvious corollary.

\begin{corollary}
\label{no-k-box}
Let $g$ be an $n$-person no-sink WTT game form and let $k$ be arbitrary, $1 \leq k \leq n$. Then $g$ contains no $k$-box.
\end{corollary}

\proof
Let us assume by contradiction that $g$ contains a $k$-box for some $1 \leq k \leq n$. Using Lemma~\ref{box-induction} we get that $g$ contains a $(k-1)$-box or a $1$-box, the latter being impossible due to Lemma~\ref{no-1-box}. Iterating the argument we subsequently get that $g$ contains a $(k-2)$-box, a $(k-3)$-box and so on, until we finally get that $g$ contains a $1$-box, which contradicts Lemma~\ref{no-1-box}.
\qed

Let us introduce more terminology connected to assignable game forms. Let $g$ be an assignable $n$-person game form and $g_i$, $1 \leq i \leq n$, functions guaranteeing the assignability of $g$. Let $x_i \in X_i$ be a strategy of player $i$. If $g_i(x_i) = a$ for some outcome $a$, we say that hyperplane $H$ defined by fixing the strategy of player $i$ to $x_i$ is {\em assigned} an outcome $a$. If a hyperplane $H$ is assigned outcome $a$ then we say that a strategy profile $x \in H$ is {\em covered} by $H$ if $g(x) = a$ and {\em not covered} by $H$ if $g(x) \neq a$.

Now let us formulate the final statement about no-sink WTT game forms.

\begin{lemma}
\label{no-sink}
Let $g$ be an $n$-person no-sink WTT game form. Then $g$ is assignable.
\end{lemma}

\proof
Let us assign to every hyperplane perpendicular to direction $i$, $1 \leq i \leq (n-1)$, its proper outcome (i.e. all hyperplanes except of those perpendicular to direction $n$ are now assigned). Let $H$ be a hyperplane perpendicular to direction $n$ and let $x,y$ be two strategy profiles in $H$ which are not covered by any hyperplane orthogonal to $H$, i.e. $g(x)$ is not the proper outcome of $H_{i}^x$ for any $1 \leq i \leq (n-1)$, and $g(y)$ is not the proper outcome of $H_{i}^y$ for any $1 \leq i \leq (n-1)$. If $g(x) \neq g(y)$ then $g$ contains a $k$-box for some $1 \leq k \leq (n-1)$ where $k$ is the number of coordinates in which $x$ and $y$ differ (they cannot differ in all $n$ coordinates since they are both in $H$). However, this is impossible in a WTT game form due to Lemma~\ref{no-k-box}, and therefore $g(x) = g(y)=a$ for some $a$ must hold. Since $x,y$ were selected as arbitrary two  not covered strategy profiles, it follows that all  not covered profiles in $H$ have the same outcome $a$ and thus can be covered by assigning $a$ to $H$. The same can be done for every hyperplane perpendicular to direction $n$ and hence $g$ is assignable.
\qed

It seems quite intuitive, that also every hyperplane perpendicular to direction $n$, which is assigned the common outcome of all uncovered profiles, is in fact assigned its proper outcome. We conjecture that it is indeed the case which would make the statement of the algorithm producing a feasible assignment much simpler: assign to each hyperplane (in any direction) its proper outcome. However, we currently have neither a proof of this fact nor a counterexample, so we leave this as an open research question. Now we are finally ready to state and prove the main result of this paper.

\begin{theorem}
\label{main-theorem}
Let $g$ be an $n$-person WTT game form. Then $g$ is assignable.
\end{theorem}

\proof
We shall proceed by induction on $n$. The base case $n=1$ is trivial. In this case $g$ is just a single line which is trivially WTT (there are no $2 \times 2$ submatrices to consider) and which is also easily assignable (each strategy of the single player is assigned the only outcome belonging to that strategy).

Let us assume for the induction step that the statement of the theorem is true for all $(n-1)$-person WTT game forms and let $g$ be an $n$-person WTT game form. Now there are two possibilities. Either $g$ has no sink hyperplane in any direction and then it is assignable by Lemma~\ref{no-sink}, or there exists a direction $i$ and a hyperplane $H_j$ perpendicular to $i$ such that $H_j$ is a sink hyperplane.  In the latter case there exists an outcome $c_k$ such that $H_k \stackrel{c_k}{\longrightarrow} H_j$ for every hyperplane $H_k$ perpendicular to direction $i$, $k \neq j$. We assign $c_k$ to $H_k$ for every $k \neq j$ which covers all profiles in $H_k^{\neq}(j)$ regions of the hyperplanes $H_k$, $k \neq j$ (the domination $H_k \stackrel{c_k}{\longrightarrow} H_j$ implies $g(x) = c_k$ for every $x \in H_k^{\neq}(j)$).

It remains to cover profiles in the $H_k^=(j)$ regions of the hyperplanes $H_k$, $k \neq j$, and in hyperplane $H_j$. However, $H_j$ is an $(n-1)$-person WTT game form which is assignable by the induction hypothesis. Moreover, if we extend the assignment of all $(n-2)$-dimensional hyperplanes inside of $H_j$ given by the hypothesis to be the same for the $(n-1)$-dimensional hyperplanes originating from the $(n-2)$-dimensional hyperplanes by adding the coordinate $i$ as a running index (extending the $(n-2)$-dimensional hyperplanes along the lines in direction $i$) then this extended assignment clearly covers all profiles in the $H_k^=(j)$ regions of all hyperplanes $H_k$, $k \neq j$. Thus all strategy profiles in $g$ are covered, which finishes the proof.
\qed

Note that the proof of Theorem~\ref{main-theorem} gives a recursive algorithm constructing a feasible assignment for an arbitrary $n$-person WTT game form. That has an impact on the complexity of recognizing the WTT property and subsequently constructing a feasible assignment for a WTT game form, as we shall see in the next section.

\section{Complexity of recognition of WTT and AS game forms}
\label{complexity-section}
First let us realize that WTT game forms can be recognized in polynomial time with respect to their size, i.e. with respect to the total number of strategy profiles. Moreover, let us note that for a WTT game form a feasible assignment can be constructed  in polynomial time as well using the algorithm implicitly present in the proof of Theorem~\ref{main-theorem}.

\begin{theorem}
\label{complexity-theorem}
Let $g$ be an $n$-person game form of size $s_1 \times s_2 \times \cdots \times s_n$. Let us denote
$s = \sum_{i=1}^n s_i$ the sum of sizes in all directions and $p = \prod_{i=1}^n s_i$ the product of sizes in all directions, i.e. let $p$ be the total number of all strategy profiles in $g$. Then it can be tested in $O(np^2)$ time whether $g$ is WTT and in the affirmative case a feasible assignment for $g$ can be constructed in $O(nsp)$ time.
\end{theorem}

\proof
To test the WTT property, it suffices to test for each direction $i$ all $2 \times 2$ subarrays defined by a choice of two of the $p/s_i$ distinct lines in direction $i$ and two of the $s_i$ distinct hyperplanes perpendicular to direction $i$. There are $O((p/s_i)^2)$ pairs of such lines and $O((s_i)^2)$ pairs of such hyperplanes. Thus there are $O(p^2)$ $2 \times 2$ subarrays to check for direction $i$ and thus altogether $O(np^2)$ subarrays for all directions (checking each $2 \times 2$ subarray takes of course just a constant time).

Now let us assume that $g$ is WTT. Given two hyperplanes perpendicular to direction $i$ (each containing $p/s_i$ strategy profiles) it takes $O(p/s_i)$ time to detect the dominance relation between them. There are $O(s_i^2)$ such pairs of hyperplanes and so it takes $O(s_i p)$ time to build the dominance graph for player $i$. Thus it takes $O(sp)$ time to build the dominance graphs for all players.

In case $g$ is a no-sink game form then, following the proof of Lemma~\ref{no-sink}, all hyperplanes perpendicular to directions other than $n$ are assigned their proper outcomes (these outcomes are contained in the dominance graphs as edge labels). Given a hyperplane $H$ perpendicular to direction $n$ and a profile $x \in H$ it takes $O(n)$ time to check whether it is covered by one of the $n-1$ already assigned hyperplanes going through $x$. Therefore it takes $O(np)$ time to check all profiles and assign outcomes also to hyperplanes perpendicular to direction $n$ which is dominated by $O(sp)$ time necessary to build the dominance graphs.

In case $g$ contains a sink hyperplane then, following the proof of Theorem~\ref{main-theorem}, a recursion is invoked. This recursion has depth at most $n$ and at each level the time needed to build the dominance graphs is of course dominated by $O(sp)$. Thus the total time needed to get a feasible assignment is $O(nsp)$.
\qed

Obviously, the above theorem is valid also for partially defined game forms. Recall that the WTT property in this case means that all undefined values are all replaced by a single extra outcome and a WTT property of this fully defined game form is required. Thus the complexity of the recognition problem is equivalent to the fully defined case. When constructing a feasible assignment for a partially defined WTT game form using the procedure for the fully defined case, some hyperplanes may be assigned the extra outcome. This just means that such hyperplanes are not needed to "cover" the defined outcomes in the partially defined game form, and hence each of these hyperplanes may be assigned an arbitrary outcome (instead of the extra outcome) in a feasible assignment of the partially defined game form.

We have seen above that WTT game forms (fully or partially defined) can be recognized in polynomial time. What is the complexity of recognition for assignable game forms? As we showed in the previous section, the set of assignable game forms is a superset of WTT game forms. In fact, it is a strict superset of the WTT ones even in the two dimensional case $n=2$. See examples in~\cite{BGMP10}, where the implication $AS \; \Rightarrow \; TT$ is disproved, or the first example in~(\ref{examples}) in the introduction. The next three subsections show how difficult it is to recognize this strict superset under different additional conditions. 

\subsection{Complexity of recognition of assignable game forms for n=2}
One way to attack the recognition problem is to formulate the assignability of a game form $g$ (both for the fully defined and partially defined cases) as a CNF satisfiability problem. If we introduce Boolean variables $y_{ij}^k$ for $1 \leq i \leq n$, $j \in X_i$, and $k \in A$, where $y_{ij}^k = 1$ means $g_i(j) = k$, then the desired CNF consists of two types of clauses. The first type guarantees for every strategy profile $x = (x_1, \ldots, x_n) \in X_1 \times \cdots \times X_n$ where $g(x) = k$ that it is "covered" by one of the separating functions:
\[ (y_{1x_1}^k \vee y_{2x_2}^k \vee \ldots \vee y_{nx_n}^k)
\]
Note that the size of each such clause (the number of literals in it) is given by the dimension of $g$ (by the number of players), and the number of such clauses is equal to the number of profiles for which $g$ is defined. The second type of clauses then guarantees that at most one outcome from $A$ is assigned to every $g_i(j)$,  $1 \leq i \leq n$, $j \in X_i$:
\[ \bigwedge_{k \neq \ell \in A}(\oy_{ij}^k \vee \oy_{ij}^{\ell})
\]
These clauses are all quadratic (two literals per clause). It is not necessary to require that exactly one outcome is assigned to every $g_i(j)$  (requiring at most one outcome suffices), because any partially defined feasible assignment can be of course arbitrarily completed to a fully defined one.

Note, that for $n=2$ the above formulation yields a 2-SAT instance (all clauses are quadratic), which immediately implies that the assignability of a two-person game form (partially or fully defined) can be recognized (and a feasible assignment constructed, if it exists) in polynomial time with respect to the size of the game form. On the other hand, given a fully defined $n$-person game form $g$ with $n \geq 3$, the complexity of recognizing whether $g$ is assignable is not known. We shall address this problem in the rest of this section. First we shall show that for partially defined game forms the recognition problem is NP-complete for $n \geq 3$. Then we will modify this proof to show that for fully defined game forms the recognition problem is NP-complete for $n \geq 4$, leaving the case $n=3$ open.

\subsection{Complexity of recognition of partially defined assignable game forms for $n\geq3$}
In this subsection we show that recognizing assignable partially defined $n$-person game forms is NP-complete already for $n=3$. This problem is obviously in NP (for any $n$) as checking a feasibility of a given assignment can be easily done in polynomial time with respect to the size of the game form. The hardness part is proved in the following theorem.

\begin{theorem}
\label{complexity-partially-defined}
It is NP-hard to recognize, whether a given partially defined 3-person game form is assignable or not.
\end{theorem}

\proof
We will proceed by constructing a reduction from the known NP-hard problem 3-SAT, a satisfiability of CNFs with exactly three literals per clause, where we also assume without any loss of generality that no clause contains two literals of the same variable. Let
\[
\Phi = \bigwedge_{i=1}^m C_i = \bigwedge_{i=1}^m (u_i \vee v_i \vee w_i)
\]
be an instance of 3-SAT, i.e. a 3-CNF on variables $x_1, \ldots ,x_n$ where each $u_i$, $v_i$, and $w_i$ is a positive or a negative occurrence of some variable. We associate outcomes $c_1, \ldots ,c_m$ with the clauses $C_1, \ldots ,C_m$ of $\Phi$ and define a partially defined 3-person game form $g_{\Phi}$. It consists of an $m \times m \times m$ box $B$, where $g_{\Phi}(i,i,i) = c_i$, i.e. box $B$ contains the outcomes $c_1, \ldots ,c_m$ on its main diagonal and it is undefined everywhere else. Let us denote $H_1^1, \ldots ,H_m^1$ the hyperplanes perpendicular to direction $1$, and similarly for directions $2$ and $3$. Box $B$ serves as a ``selector''. The fact that hyperplane $H_i^1$ is assigned outcome $c_i$ means that clause $C_i$ is satisfied by literal $u_i$ (the literal $u_i$ gets value true), and similarly for $H_i^2$ and $v_i$ and also $H_i^3$ and $w_i$. Clearly, $B$ by itself is assignable in many ways - each strategy profile on the main diagonal of $B$ can be covered by any of the three hyperplanes incident with it. However, not every such assignment corresponds to a truth assignment to variables $x_1, \ldots ,x_n$ as it disregards the fact that two distinct literals may share the same variable. To establish a one-to-one correspondence between feasible assignments of $g_{\Phi}$ and satisfying truth assignments of $\Phi$ we will add ``gadgets'' to the box $B$ which will guarantee that:
\begin{enumerate}
\item If $u_i = u_j$ for $i \neq j$, i.e. both literals are two occurrences of the same variable with the same polarity, then either $H_i^1$ is assigned $c_i$ and simultaneously $H_j^1$ is assigned $c_j$ (both $u_i$ and $u_j$ are true), or neither $H_i^1$ is assigned $c_i$ nor $H_j^1$ is assigned $c_j$ (both $u_i$ and $u_j$ are false). Similarly for $v_i = v_j$ and $w_i = w_j$. In each of these three cases we need to force the assigned outcomes in the above described way for a pair of parallel hyperplanes.
\item If $u_i = v_j$ for $i \neq j$, i.e. both literals are two occurrences of the same variable with the same polarity, then either $H_i^1$ is assigned $c_i$ and simultaneously $H_j^2$ is assigned $c_j$ (both $u_i$ and $u_j$ are true), or neither $H_i^1$ is assigned $c_i$ nor $H_j^2$ is assigned $c_j$ (both $u_i$ and $u_j$ are false). Similarly for $u_i = w_j$ and $v_i = w_j$. In each of these three cases we need to force the assigned outcomes in the above described way for a pair of perpendicular hyperplanes.
\item If $u_i = \ou_j$ for $i \neq j$, i.e. both literals are two occurrences of the same variable with complementary polarities, then either $H_i^1$ is assigned $c_i$ or $H_j^1$ is assigned $c_j$ but not both (exactly one of $u_i$ and $u_j$ is true and exactly one false). Similarly for $v_i = \ov_j$ and $w_i = \ow_j$. In each of these three cases we need to force the assigned outcomes in the above described way for a pair of parallel hyperplanes.
\item If $u_i = \ov_j$ for $i \neq j$, i.e. both literals are two occurrences of the same variable with complementary polarities, then either $H_i^1$ is assigned $c_i$ or $H_j^2$ is assigned $c_j$ but not both (exactly one of $u_i$ and $v_j$ is true and exactly one false). Similarly for $u_i = \ow_j$ and $v_i = \ow_j$. In each of these three cases we need to force the assigned outcomes in the above described way for a pair of perpendicular hyperplanes.
\end{enumerate}
Each such gadget will add a constant number of ``dedicated'' hyperplanes outside of box $B$, where ``dedicated'' means that no two gadgets share a common added hyperplane (of course they may share hyperplanes incident to box $B$). Since the number of the above pair-wise relations is at most quadratic in the size of $\Phi$, the constructed game form $g_{\Phi}$ has a polynomial size with respect to the size of $\Phi$.

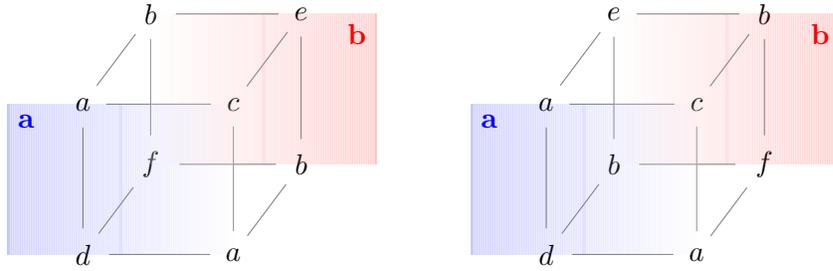
\begin{figure}[ht]
\centering
\begin{tikzpicture}[z={(0.45cm,0.6cm)},x={(1cm,0cm)},y={(0cm,1cm)},scale=2]

\foreach \x in {1,2} \foreach \y in {1,2}
	\draw[gray,thin] (\x,\y,1) -- (\x,\y,2);

\foreach \x in {1,2} \foreach \z in {1,2}
	\draw[gray,thin] (\x,1,\z) -- (\x,2,\z);
	
\foreach \y in {1,2} \foreach \z in {1,2}
	\draw[gray,thin] (1,\y,\z) -- (2,\y,\z);
	


\node[circle,fill=white] at (1,1,1) {$d$};
\node[circle,fill=white] at (1,2,1) {$a$};
\node[circle,fill=white] at (2,1,1) {$a$};
\node[circle,fill=white] at (2,2,1) {$c$};
\node[circle,fill=white] at (1,1,2) {$f$};
\node[circle,fill=white] at (1,2,2) {$b$};
\node[circle,fill=white] at (2,1,2) {$b$};
\node[circle,fill=white] at (2,2,2) {$e$};

\shade[left color=blue,right color=white,opacity=0.1] (0.5,2,1) -- (2,2,1) -- (2,1,1) -- (0.5,1,1) -- cycle;

\node[below right] at (0.5,2,1) {$\color{blue}\mathbf{a}$};

\shade[left color=white,right color=red,opacity=0.1] (1,1,2) -- (1,2,2) -- (2.5,2,2) -- (2.5,1,2) -- cycle;

\node[below left] at (2.5,2,2) {$\color{red}\mathbf{b}$};

\end{tikzpicture}
\hspace*{10mm}
\begin{tikzpicture}[z={(0.45cm,0.6cm)},x={(1cm,0cm)},y={(0cm,1cm)},scale=2]

\foreach \x in {1,2} \foreach \y in {1,2}
	\draw[gray,thin] (\x,\y,1) -- (\x,\y,2);

\foreach \x in {1,2} \foreach \z in {1,2}
	\draw[gray,thin] (\x,1,\z) -- (\x,2,\z);
	
\foreach \y in {1,2} \foreach \z in {1,2}
	\draw[gray,thin] (1,\y,\z) -- (2,\y,\z);
	


\node[circle,fill=white] at (1,1,1) {$d$};
\node[circle,fill=white] at (1,2,1) {$a$};
\node[circle,fill=white] at (2,1,1) {$a$};
\node[circle,fill=white] at (2,2,1) {$c$};
\node[circle,fill=white] at (1,1,2) {$b$};
\node[circle,fill=white] at (1,2,2) {$e$};
\node[circle,fill=white] at (2,1,2) {$f$};
\node[circle,fill=white] at (2,2,2) {$b$};

\shade[left color=blue,right color=white,opacity=0.1] (0.5,2,1) -- (2,2,1) -- (2,1,1) -- (0.5,1,1) -- cycle;

\node[below right] at (0.5,2,1) {$\color{blue}\mathbf{a}$};

\shade[left color=white,right color=red,opacity=0.1] (1,1,2) -- (1,2,2) -- (2.5,2,2) -- (2.5,1,2) -- cycle;

\node[below left] at (2.5,2,2) {$\color{red}\mathbf{b}$};

\end{tikzpicture}
\caption{Two forcing cubes: Here we assume that $a$, $b$, $c$ $d$, $e$, and $f$ are six distinct outcomes. Since they are covered by six planes, both copies of $a$ must be covered by the front plane, and both copies of $b$ must be covered by the back plane. Consequently, we must have $a$ assigned to the front and $b$ to the back planes in all feasible assignments. \label{f1}}
\end{figure}

Before constructing the four types of gadgets with above described properties, let us construct common ``forcing'' components of such gadgets (called {\em forcing cubes}). These are $2 \times 2 \times 2$ arrays with six distinct outcomes $a,b,c,d,e,f$ arranged in the eight corners (strategy profiles) of the array in one of the two possible ways shown in Figure \ref{f1}.

In both forcing cubes there are six hyperplanes to cover all eight corners with six distinct outcomes. Thus, both $a$ outcomes must be covered by the same hyperplane and so do both $b$ outcomes. Hence, in both forcing cubes the front hyperplane is forced to be assigned $a$ and the back hyperplane is forced to be assigned $b$ in any feasible assignment.

Let us now construct the four types of required gadgets.
\begin{enumerate}
\item Let $u_i = u_j$. Let us add two hyperplanes $H_{\ell}^2$ and $H_{\ell+1}^2$ perpendicular to direction $2$ and two hyperplanes $H_k^3$ and $H_{k+1}^3$ perpendicular to direction $3$ for some $k,\ell > m$. Let us consider four distinct outcomes $a,b,d,e$ not contained among $c_1, \ldots ,c_m$, and the $2 \times 2 \times 2$ array defined by the intersections of the four added hyperplanes with $H_i^1$ and $H_i^2$ as in Figure \ref{f-gadget1}.

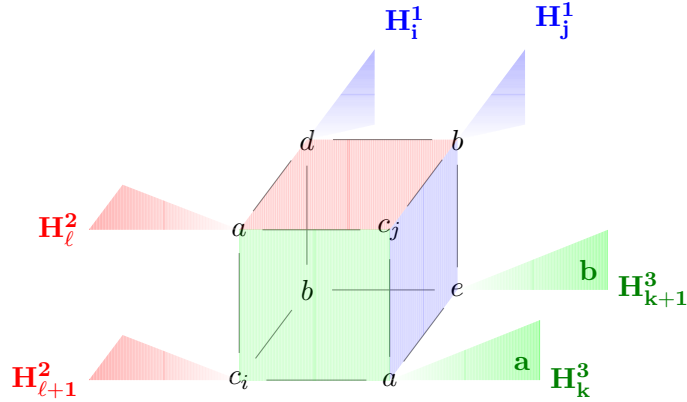
\begin{figure}[ht]
\centering
\begin{tikzpicture}[z={(0.45cm,0.6cm)},x={(1cm,0cm)},y={(0cm,1cm)},scale=2]

\foreach \x in {1,2} \foreach \y in {1,2}
	\draw[gray,thin] (\x,\y,1) -- (\x,\y,2);

\foreach \x in {1,2} \foreach \z in {1,2}
	\draw[gray,thin] (\x,1,\z) -- (\x,2,\z);
	
\foreach \y in {1,2} \foreach \z in {1,2}
	\draw[gray,thin] (1,\y,\z) -- (2,\y,\z);
	




\node[circle,fill=white] at (1,1,1) {$c_i$};
\node[circle,fill=white] at (1,2,1) {$a$};
\node[circle,fill=white] at (2,1,1) {$a$};
\node[circle,fill=white] at (2,2,1) {$c_j$};
\node[circle,fill=white] at (1,1,2) {$b$};
\node[circle,fill=white] at (1,2,2) {$d$};
\node[circle,fill=white] at (2,1,2) {$e$};
\node[circle,fill=white] at (2,2,2) {$b$};

\shade[right color=blue,left color=white,top color=blue,bottom color=white,opacity=0.2] (1,2,2) -- (1,2,3) -- (1,1.5,3) -- cycle;

\node[above right] at (1,2,3) {$\color{blue}\mathbf{H^1_i}$};

\shade[right color=blue,left color=white,top color=blue,bottom color=white,opacity=0.2] (2,2,2) -- (2,2,3) -- (2,1.5,3) -- cycle;

\node[above right] at (2,2,3) {$\color{blue}\mathbf{H^1_j}$};

\node[left] at (0,2,1) {$\color{red}\mathbf{H^2_\ell}$};

\node[left] at (0,1,1) {$\color{red}\mathbf{H^2_{\ell+1}}$};

\shade[left color=red, bottom color=red, right color=white, opacity=0.2] (1,1,1) -- (0,1,1) -- (0,1,1.5) -- cycle;

\shade[left color=red, bottom color=red, right color=white, opacity=0.2] (1,2,1) -- (0,2,1) -- (0,2,1.5) -- cycle;

\node[right] at (3,1,2) {$\color{green!50!black}\mathbf{H^3_{k+1}}$};

\node[right] at (3,1,1) {$\color{green!50!black}\mathbf{H^3_k}$};

\shade[bottom color=green, right color=green, left color=white, opacity=0.2] (2,1,2) -- (3,1,2) -- (3,1.4,2) -- cycle;

\shade[bottom color=green, right color=green, left color=white, opacity=0.2] (2,1,1) -- (3,1,1) -- (3,1.4,1) -- cycle;

\shade[left color=green,right color=green,opacity=0.1] (1,1,1) -- (1,2,1) -- (2,2,1) -- (2,1,1) -- cycle;

\shade[left color=blue,right color=blue,opacity=0.1] (2,1,1) -- (2,2,1) -- (2,2,2) -- (2,1,2) -- cycle;

\shade[left color=red,right color=red,opacity=0.1] (1,2,1) -- (2,2,1) -- (2,2,2) -- (1,2,2) -- cycle;

\node[above left] at (3,1,2) {$\color{green!50!black}\mathbf{b}$};
\node[above left] at (3,1,1) {$\color{green!50!black}\mathbf{a}$};

\end{tikzpicture}
\caption{Gadget 1: This is an instance of one of the forcing cubes, thus assignments $a\to H^3_k$ and $b\to H^3_{k+1}$ are implied in any feasible assignment. These leave only two possible cyclic feasible assignments for the remaining four planes. Either $c_i\to H^1_i$, $d\to H^2_\ell$, $c_j\to H^1_j$ and $e\to H^2_{\ell+1}$, or
$d\to H^1_i$, $c_j\to H^2_\ell$, $e\to H^1_j$ and $c_i\to H^2_{\ell+1}$. Consequently, in all feasible assignments we have either both $c_i\to H^1_i$ and $c_j\to H^1_j$ or we have neither one, simultaneously.\label{f-gadget1}}
\end{figure}
This gadget is a forcing cube which forces $H_k^3$ to be assigned $a$ and $H_{k+1}^3$ to be assigned $b$. Now there are only two ways how the remaining four hyperplanes can cover the four distinct outcomes $c_i,c_j,d,e$. If $H_i^1$ is assigned $c_i$, it forces $H_{\ell}^2$ to be assigned $d$, which in turn forces $H_j^1$ to be assigned $c_j$, which finally forces $H_{\ell+1}^2$ to be assigned $e$. On the other hand, if $H_i^1$ is assigned $d$, it forces $H_{\ell+1}^2$ to be assigned $c_i$, which in turn forces $H_j^1$ to be assigned $e$, which finally forces $H_{\ell}^2$ to be assigned $c_j$. Thus the constructed gadget fulfills exactly the required properties. Note also, that the hyperplanes $H_i^1, H_j^1$ which are incident to the selector box $B$ are forced to be assigned one of the outcomes $c_i, c_j, d, e$ in every feasible assignment.
\item Let $u_i = v_j$. Let us add one hyperplane $H_{\ell}^1$ perpendicular to direction $1$, one hyperplane $H_p^2$ perpendicular to direction $2$, and two hyperplanes $H_k^3$ and $H_{k+1}^3$ perpendicular to direction $3$ for some $k,\ell,p > m$. Let us consider four distinct outcomes $a,b,d,e$ not contained among $c_1, \ldots ,c_m$ and the $2 \times 2 \times 2$ array defined by the intersections of the four added hyperplanes with $H_i^1$ and $H_j^2$. We assign the outcomes to this $2 \times 2 \times 2$ array as in Figure \ref{f-gadget2}.
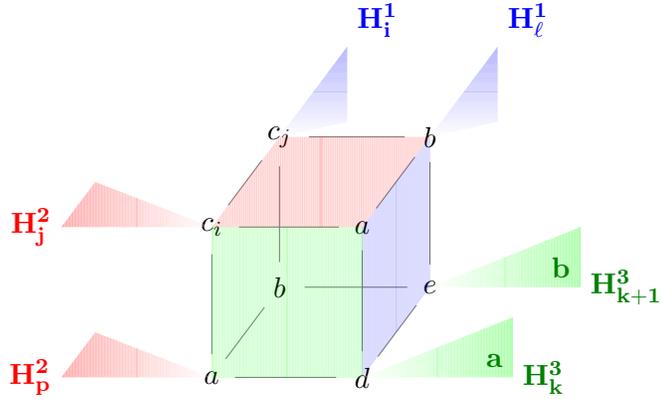
\begin{figure}[ht]
\centering
\begin{tikzpicture}[z={(0.45cm,0.6cm)},x={(1cm,0cm)},y={(0cm,1cm)},scale=2]

\foreach \x in {1,2} \foreach \y in {1,2}
	\draw[gray,thin] (\x,\y,1) -- (\x,\y,2);

\foreach \x in {1,2} \foreach \z in {1,2}
	\draw[gray,thin] (\x,1,\z) -- (\x,2,\z);
	
\foreach \y in {1,2} \foreach \z in {1,2}
	\draw[gray,thin] (1,\y,\z) -- (2,\y,\z);
	




\node[circle,fill=white] at (1,1,1) {$a$};
\node[circle,fill=white] at (1,2,1) {$c_i$};
\node[circle,fill=white] at (2,1,1) {$d$};
\node[circle,fill=white] at (2,2,1) {$a$};
\node[circle,fill=white] at (1,1,2) {$b$};
\node[circle,fill=white] at (1,2,2) {$c_j$};
\node[circle,fill=white] at (2,1,2) {$e$};
\node[circle,fill=white] at (2,2,2) {$b$};

\shade[right color=blue,left color=white,top color=blue,bottom color=white,opacity=0.2] (1,2,2) -- (1,2,3) -- (1,1.5,3) -- cycle;

\node[above right] at (1,2,3) {$\color{blue}\mathbf{H^1_i}$};

\shade[right color=blue,left color=white,top color=blue,bottom color=white,opacity=0.2] (2,2,2) -- (2,2,3) -- (2,1.5,3) -- cycle;

\node[above right] at (2,2,3) {$\color{blue}\mathbf{H^1_\ell}$};

\node[left] at (0,2,1) {$\color{red}\mathbf{H^2_j}$};

\node[left] at (0,1,1) {$\color{red}\mathbf{H^2_p}$};

\shade[left color=red, bottom color=red, right color=white, opacity=0.2] (1,1,1) -- (0,1,1) -- (0,1,1.5) -- cycle;

\shade[left color=red, bottom color=red, right color=white, opacity=0.2] (1,2,1) -- (0,2,1) -- (0,2,1.5) -- cycle;

\node[right] at (3,1,2) {$\color{green!50!black}\mathbf{H^3_{k+1}}$};

\node[right] at (3,1,1) {$\color{green!50!black}\mathbf{H^3_k}$};

\shade[bottom color=green, right color=green, left color=white, opacity=0.2] (2,1,2) -- (3,1,2) -- (3,1.4,2) -- cycle;

\shade[bottom color=green, right color=green, left color=white, opacity=0.2] (2,1,1) -- (3,1,1) -- (3,1.4,1) -- cycle;

\shade[left color=green,right color=green,opacity=0.1] (1,1,1) -- (1,2,1) -- (2,2,1) -- (2,1,1) -- cycle;

\shade[left color=blue,right color=blue,opacity=0.1] (2,1,1) -- (2,2,1) -- (2,2,2) -- (2,1,2) -- cycle;

\shade[left color=red,right color=red,opacity=0.1] (1,2,1) -- (2,2,1) -- (2,2,2) -- (1,2,2) -- cycle;

\node[above left] at (3,1,2) {$\color{green!50!black}\mathbf{b}$};
\node[above left] at (3,1,1) {$\color{green!50!black}\mathbf{a}$};

\end{tikzpicture}
\caption{Gadget 2: This is an instance of one of the forcing cubes, thus assignments $a\to H^3_k$ and $b\to H^3_{k+1}$ are implied in any feasible assignment. These leave only two possible feasible assignments for $H^1_i$ and $H^2_j$. We have either $c_i\to H^1_i$ and $c_j\to H^2_j$, or
$c_j\to H^1_i$ and $c_i\to H^2_j$. Consequently, in all feasible assignments we have either both $c_i\to H^1_i$ and $c_j\to H^2_j$ or we have neither one, simultaneously.\label{f-gadget2}}
\end{figure}

This is again a forcing cube which forces $H_k^3$ to be assigned $a$ and $H_{k+1}^3$ to be assigned $b$. Clearly, either $H_i^1$ is assigned $c_i$ and $H_j^2$ is assigned $c_j$ or alternatively $H_i^1$ is assigned $c_j$ and $H_j^2$ is assigned $c_i$ (and the outcomes $d$ and $e$ are similarly covered by $H_{\ell}^1$ and $H_p^2$ in one of the two possible ways). Thus the constructed gadget fulfills exactly the required properties. Note also, that the hyperplanes $H_i^1, H_j^2$ which are incident to the selector box $B$ are forced to be assigned one of the outcomes $c_i, c_j$ in every feasible assignment.
\item Let $u_i = \ou_j$. Let us add two hyperplanes $H_{\ell}^1$ and $H_{\ell+1}^1$ perpendicular to direction $1$, two hyperplanes $H_p^2$ and $H_{p+1}^2$ perpendicular to direction $2$, and two hyperplanes $H_k^3$ and $H_{k+1}^3$ perpendicular to direction $3$ for some $k,\ell,p > m$. Let us assign six distinct outcomes $a,b,c,d,e,f$ not contained among $c_1, \ldots ,c_m$ to the
$2 \times 2 \times 2$ array defined by the intersections of the six added hyperplanes in such a way that we get a forcing cube which forces $H_k^3$ to be assigned outcome $a$. Let us add two more hyperplanes $H_{p+2}^2$ and $H_{p+3}^2$ perpendicular to direction $2$ and consider the intersections of $H_k^3$ with $H_i^1$, $H_j^1$, $H_{p+2}^2$, and $H_{p+3}^2$. Let us assign outcomes to this $2 \times 2$ subarray as in Figure \ref{f-gadget3}.
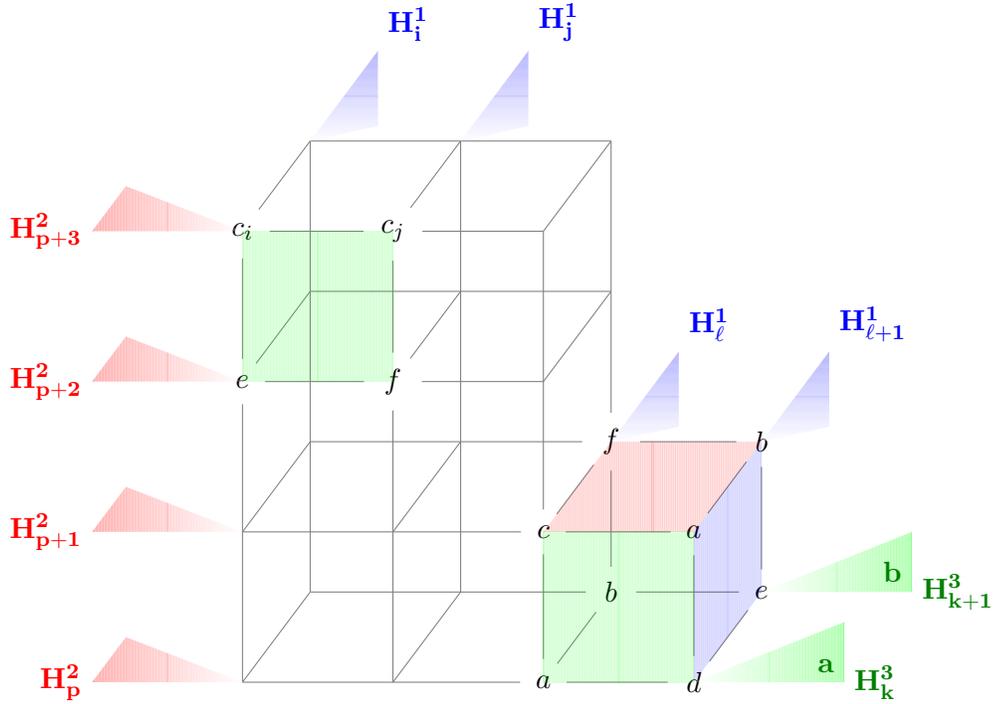
\begin{figure}[ht]
\centering
\begin{tikzpicture}[z={(0.45cm,0.6cm)},x={(1cm,0cm)},y={(0cm,1cm)},scale=2]

\foreach \x in {-1,0,1} \foreach \y in {1,2,3,4}
	\draw[gray,thin] (\x,\y,1) -- (\x,\y,2);

\foreach \x in {-1,0,1} \foreach \z in {1,2}
	\draw[gray,thin] (\x,1,\z) -- (\x,4,\z);
	
\foreach \y in {1,2,3,4} \foreach \z in {1,2}
	\draw[gray,thin] (-1,\y,\z) -- (1,\y,\z);
	
\foreach \y in {1,2} \foreach \z in {1,2}
	\draw[gray,thin] (1,\y,\z) -- (2,\y,\z);

\draw[gray,thin] (2,1,1) -- (2,2,1) -- (2,2,2) -- (2,1,2) -- cycle;

\node[circle,fill=white] at (1,1,1) {$a$};
\node[circle,fill=white] at (1,2,1) {$c$};
\node[circle,fill=white] at (2,1,1) {$d$};
\node[circle,fill=white] at (2,2,1) {$a$};
\node[circle,fill=white] at (1,1,2) {$b$};
\node[circle,fill=white] at (1,2,2) {$f$};
\node[circle,fill=white] at (2,1,2) {$e$};
\node[circle,fill=white] at (2,2,2) {$b$};

\node[circle,fill=white] at (-1,3,1) {$e$};
\node[circle,fill=white] at (-1,4,1) {$c_i$};
\node[circle,fill=white] at (0,3,1) {$f$};
\node[circle,fill=white] at (0,4,1) {$c_j$};

\shade[right color=blue,left color=white,top color=blue,bottom color=white,opacity=0.2] (1,2,2) -- (1,2,3) -- (1,1.5,3) -- cycle;

\node[above right] at (1,2,3) {$\color{blue}\mathbf{H^1_\ell}$};

\shade[right color=blue,left color=white,top color=blue,bottom color=white,opacity=0.2] (2,2,2) -- (2,2,3) -- (2,1.5,3) -- cycle;

\node[above right] at (2,2,3) {$\color{blue}\mathbf{H^1_{\ell+1}}$};

\shade[right color=blue,left color=white,top color=blue,bottom color=white,opacity=0.2] (-1,4,2) -- (-1,4,3) -- (-1,3.5,3) -- cycle;

\node[above right] at (-1,4,3) {$\color{blue}\mathbf{H^1_i}$};

\shade[right color=blue,left color=white,top color=blue,bottom color=white,opacity=0.2] (0,4,2) -- (0,4,3) -- (0,3.5,3) -- cycle;

\node[above right] at (0,4,3) {$\color{blue}\mathbf{H^1_j}$};

\node[left] at (-2,4,1) {$\color{red}\mathbf{H^2_{p+3}}$};

\node[left] at (-2,3,1) {$\color{red}\mathbf{H^2_{p+2}}$};

\node[left] at (-2,2,1) {$\color{red}\mathbf{H^2_{p+1}}$};

\node[left] at (-2,1,1) {$\color{red}\mathbf{H^2_p}$};

\shade[left color=red, bottom color=red, right color=white, opacity=0.2] (-1,1,1) -- (-2,1,1) -- (-2,1,1.5) -- cycle;

\shade[left color=red, bottom color=red, right color=white, opacity=0.2] (-1,2,1) -- (-2,2,1) -- (-2,2,1.5) -- cycle;

\shade[left color=red, bottom color=red, right color=white, opacity=0.2] (-1,3,1) -- (-2,3,1) -- (-2,3,1.5) -- cycle;

\shade[left color=red, bottom color=red, right color=white, opacity=0.2] (-1,4,1) -- (-2,4,1) -- (-2,4,1.5) -- cycle;

\node[right] at (3,1,2) {$\color{green!50!black}\mathbf{H^3_{k+1}}$};

\node[right] at (3,1,1) {$\color{green!50!black}\mathbf{H^3_k}$};

\shade[bottom color=green, right color=green, left color=white, opacity=0.2] (2,1,2) -- (3,1,2) -- (3,1.4,2) -- cycle;

\shade[bottom color=green, right color=green, left color=white, opacity=0.2] (2,1,1) -- (3,1,1) -- (3,1.4,1) -- cycle;

\shade[left color=green,right color=green,opacity=0.1] (1,1,1) -- (1,2,1) -- (2,2,1) -- (2,1,1) -- cycle;

\shade[left color=blue,right color=blue,opacity=0.1] (2,1,1) -- (2,2,1) -- (2,2,2) -- (2,1,2) -- cycle;

\shade[left color=red,right color=red,opacity=0.1] (1,2,1) -- (2,2,1) -- (2,2,2) -- (1,2,2) -- cycle;

\shade[left color=green,right color=green,opacity=0.1] (-1,3,1) -- (-1,4,1) -- (0,4,1) -- (0,3,1) -- cycle;

\node[above left] at (3,1,2) {$\color{green!50!black}\mathbf{b}$};
\node[above left] at (3,1,1) {$\color{green!50!black}\mathbf{a}$};

\end{tikzpicture}
\caption{Gadget 3: In this gadget we have a forcing cube disjoint from the hyperplanes incident with the selector box, implying that $a$ is assigned to $H^3_k$ in all feasible assignments. The intersection of this hyperplane with four others, as in the picture above, reduces the possible feasible assignments to the cases when either $c_i$ is assigned to $H^1_i$ or $c_j$ is assigned to $H^1_j$, but not both.\label{f-gadget3}}
\end{figure}

Since $H_k^3$ is forced to be assigned outcome $a$, there are just two ways how to cover the four outcomes in this $2 \times 2$ subarray. Either $H_i^1$ is assigned $c_i$, $H_{p+3}^2$ is assigned $e$, $H_j^1$ is assigned $f$, and $H_{p+2}^2$ is assigned $c_j$, or alternatively $H_i^1$ is assigned $e$, $H_{p+2}^2$ is assigned $c_i$, $H_j^1$ is assigned $c_j$, and $H_{p+3}^2$ is assigned $f$. Thus the constructed gadget fulfills exactly the required properties. Note also, that the hyperplanes $H_i^1, H_j^1$ which are incident to the selector box $B$ are forced to be assigned one of the outcomes $c_i, c_j, e, f$ in every feasible assignment.
\item Let $u_i = \ov_j$. Let us add two hyperplanes $H_{\ell}^1$ and $H_{\ell+1}^1$ perpendicular to direction $1$, two hyperplanes $H_p^2$ and $H_{p+1}^2$ perpendicular to direction $2$, and two hyperplanes $H_k^3$ and $H_{k+1}^3$ perpendicular to direction $3$ for some $k,\ell,p > m$. Let us assign six distinct outcomes $a,b,c,d,e,f$ not contained among $c_1, \ldots ,c_m$ to the
$2 \times 2 \times 2$ array defined by the intersections of the six added hyperplanes in such a way, that we get a forcing cube which forces $H_k^3$ to be assigned outcome $a$. Let us add one more hyperplane $H_{\ell+2}^1$ perpendicular to direction $1$ and one more hyperplane $H_{p+2}^2$  perpendicular to direction $2$ and consider the intersections of $H_k^3$ with $H_i^1$, $H_j^2$, $H_{\ell+2}^1$, and $H_{p+2}^2$. Let us assign outcomes to this $2 \times 2$ subarray as in Figure \ref{f-gadget4}.
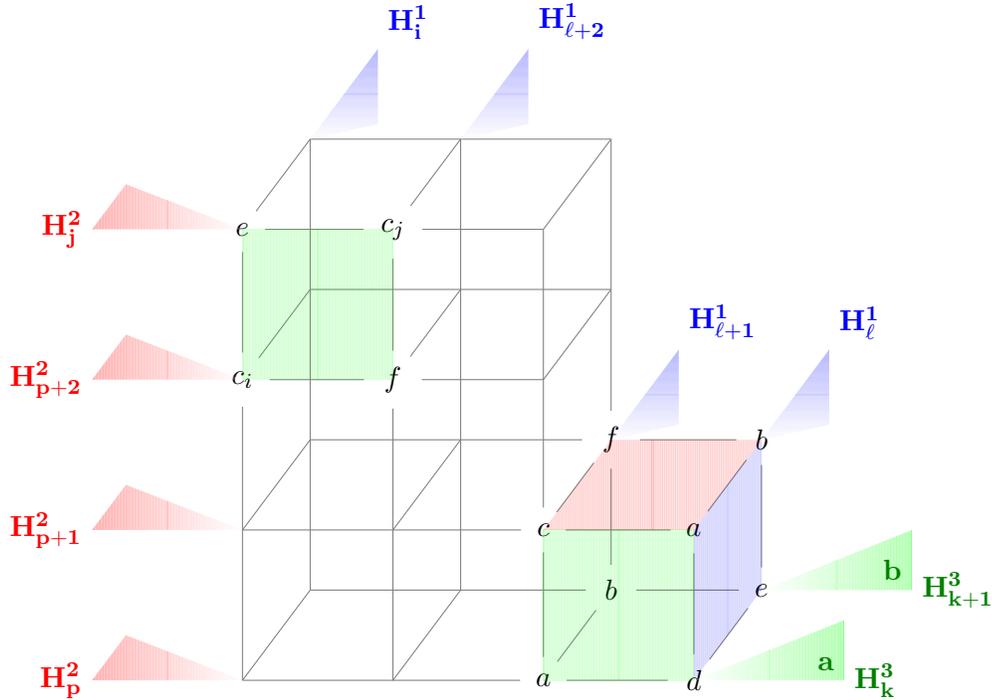
\begin{figure}[ht]
\centering
\begin{tikzpicture}[z={(0.45cm,0.6cm)},x={(1cm,0cm)},y={(0cm,1cm)},scale=2]

\foreach \x in {-1,0,1} \foreach \y in {1,2,3,4}
	\draw[gray,thin] (\x,\y,1) -- (\x,\y,2);

\foreach \x in {-1,0,1} \foreach \z in {1,2}
	\draw[gray,thin] (\x,1,\z) -- (\x,4,\z);
	
\foreach \y in {1,2,3,4} \foreach \z in {1,2}
	\draw[gray,thin] (-1,\y,\z) -- (1,\y,\z);
	
\foreach \y in {1,2} \foreach \z in {1,2}
	\draw[gray,thin] (1,\y,\z) -- (2,\y,\z);

\draw[gray,thin] (2,1,1) -- (2,2,1) -- (2,2,2) -- (2,1,2) -- cycle;

\node[circle,fill=white] at (1,1,1) {$a$};
\node[circle,fill=white] at (1,2,1) {$c$};
\node[circle,fill=white] at (2,1,1) {$d$};
\node[circle,fill=white] at (2,2,1) {$a$};
\node[circle,fill=white] at (1,1,2) {$b$};
\node[circle,fill=white] at (1,2,2) {$f$};
\node[circle,fill=white] at (2,1,2) {$e$};
\node[circle,fill=white] at (2,2,2) {$b$};

\node[circle,fill=white] at (-1,3,1) {$c_i$};
\node[circle,fill=white] at (-1,4,1) {$e$};
\node[circle,fill=white] at (0,3,1) {$f$};
\node[circle,fill=white] at (0,4,1) {$c_j$};

\shade[right color=blue,left color=white,top color=blue,bottom color=white,opacity=0.2] (1,2,2) -- (1,2,3) -- (1,1.5,3) -- cycle;

\node[above right] at (1,2,3) {$\color{blue}\mathbf{H^1_{\ell+1}}$};

\shade[right color=blue,left color=white,top color=blue,bottom color=white,opacity=0.2] (2,2,2) -- (2,2,3) -- (2,1.5,3) -- cycle;

\node[above right] at (2,2,3) {$\color{blue}\mathbf{H^1_{\ell}}$};

\shade[right color=blue,left color=white,top color=blue,bottom color=white,opacity=0.2] (-1,4,2) -- (-1,4,3) -- (-1,3.5,3) -- cycle;

\node[above right] at (-1,4,3) {$\color{blue}\mathbf{H^1_i}$};

\shade[right color=blue,left color=white,top color=blue,bottom color=white,opacity=0.2] (0,4,2) -- (0,4,3) -- (0,3.5,3) -- cycle;

\node[above right] at (0,4,3) {$\color{blue}\mathbf{H^1_{\ell+2}}$};

\node[left] at (-2,4,1) {$\color{red}\mathbf{H^2_j}$};

\node[left] at (-2,3,1) {$\color{red}\mathbf{H^2_{p+2}}$};

\node[left] at (-2,2,1) {$\color{red}\mathbf{H^2_{p+1}}$};

\node[left] at (-2,1,1) {$\color{red}\mathbf{H^2_p}$};

\shade[left color=red, bottom color=red, right color=white, opacity=0.2] (-1,1,1) -- (-2,1,1) -- (-2,1,1.5) -- cycle;

\shade[left color=red, bottom color=red, right color=white, opacity=0.2] (-1,2,1) -- (-2,2,1) -- (-2,2,1.5) -- cycle;

\shade[left color=red, bottom color=red, right color=white, opacity=0.2] (-1,3,1) -- (-2,3,1) -- (-2,3,1.5) -- cycle;

\shade[left color=red, bottom color=red, right color=white, opacity=0.2] (-1,4,1) -- (-2,4,1) -- (-2,4,1.5) -- cycle;

\node[right] at (3,1,2) {$\color{green!50!black}\mathbf{H^3_{k+1}}$};

\node[right] at (3,1,1) {$\color{green!50!black}\mathbf{H^3_k}$};

\shade[bottom color=green, right color=green, left color=white, opacity=0.2] (2,1,2) -- (3,1,2) -- (3,1.4,2) -- cycle;

\shade[bottom color=green, right color=green, left color=white, opacity=0.2] (2,1,1) -- (3,1,1) -- (3,1.4,1) -- cycle;

\shade[left color=green,right color=green,opacity=0.1] (1,1,1) -- (1,2,1) -- (2,2,1) -- (2,1,1) -- cycle;

\shade[left color=blue,right color=blue,opacity=0.1] (2,1,1) -- (2,2,1) -- (2,2,2) -- (2,1,2) -- cycle;

\shade[left color=red,right color=red,opacity=0.1] (1,2,1) -- (2,2,1) -- (2,2,2) -- (1,2,2) -- cycle;

\shade[left color=green,right color=green,opacity=0.1] (-1,3,1) -- (-1,4,1) -- (0,4,1) -- (0,3,1) -- cycle;

\node[above left] at (3,1,2) {$\color{green!50!black}\mathbf{b}$};
\node[above left] at (3,1,1) {$\color{green!50!black}\mathbf{a}$};

\end{tikzpicture}
\caption{Gadget 4: In this gadget we have a forcing cube disjoint from the hyperplanes incident with the selector box, implying that $a$ is assigned to $H^3_k$ in all feasible assignments. The intersection of this hyperplane with four others, as in the picture above, reduces the possible feasible assignments to the cases when either $c_i$ is assigned to $H^1_i$ or $c_j$ is assigned to $H^1_j$, but not both.\label{f-gadget4}}
\end{figure}

Since $H_k^3$ is forced to be assigned outcome $a$, there are just two ways how to cover the four outcomes in this $2 \times 2$ subarray. Either $H_i^1$ is assigned $c_i$, $H_j^2$ is assigned $e$, $H_{\ell+2}^1$ is assigned $c_j$, and $H_{p+2}^2$ is assigned $f$, or alternatively $H_i^1$ is assigned $e$, $H_{p+2}^2$ is assigned $c_i$, $H_{\ell+2}^1$ is assigned $f$, and $H_j^2$ is assigned $c_j$. Thus the constructed gadget fulfills exactly the required properties. Note also, that the hyperplanes $H_i^1, H_j^2$ which are incident to the selector box $B$ are forced to be assigned one of the outcomes $c_i, c_j, e, f$ in every feasible assignment.
\end{enumerate}

It follows from the above constructions that if there exists a feasible assignment to game form $g_\Phi$ then $\Phi$ has a satisfying assignment. Indeed, the feasibility of the assignment for the strategy profiles on the main diagonal of the selector box implies that each clause has a satisfying literal, and the feasibility of the assignment in the strategy profiles of the gadgets imply that these truth values are consistent.

Conversely, if we have a satisfying truth assignment to $\Phi$, then we can derive a feasible assignment to all hyperplanes $H^1_i$, $H^2_j$ and $H^3_k$ which cover the strategy profiles along the diagonal of the selector box, and extend these to cover all other strategy profiles by the proven properties of the gadgets.
\qed

\subsection{Complexity of recognition of fully defined assignable game forms for $n\geq4$}
In this subsection we will modify the proof of Theorem~\ref{complexity-partially-defined} to show that the recognition problem is NP-complete also for fully defined game forms, this time for $n \geq 4$, leaving the case $n=3$ open.

\begin{theorem}
\label{complexity-fully-defined}
It is NP-hard to recognize, whether a given fully defined 4-person game form is assignable or not.
\end{theorem}

\proof
Let us repeat the construction from the proof of Theorem~\ref{complexity-partially-defined} with these changes:
\begin{itemize}
\item All strategy profiles which were undefined in the construction now get a new additional outcome $*$, which produces a fully defined game form.
\item If the construction produced a 3-person game form of size $s_1 \times s_2 \times s_3$ we shall consider it now as a 4-person game form of size $s_1 \times s_2 \times s_3 \times 1$, and denote the single hyperplane perpendicular to the added direction $H_1^4$.
\item We shall assume that the input 3-CNF $\Phi$ satisfies the following additional property: if we delete any two clauses from $\Phi$, the remaining 3-CNF contains some clause $C_i$ with a non-trivial literal in the first position, clause $C_j$ with a non-trivial literal in the second position, and clause $C_k$ with a non-trivial literal in the third position, where a trivial literal is a literal which represents the only occurrence of its variable in the entire formula. This assumption can be made without losing the NP-hardness of the 3-SAT problem restricted to such inputs. Let us first note that we can assume that every variable appears at least twice in two different clauses of the input 3-CNF. Otherwise, we can fix the value of the unique appearance without changing satisfiability of the input.
We claim that a 3-CNF that does not satisfy the property claimed above, and that has every variable appearing at least twice, cannot have more than $10$ clauses. To see this let us assume that by deleting the first two clauses we have a trivial literal in each of the remaining clauses. All of the variables of these trivial literals must have then their second appearance in the deleted clauses, that is we cannot have more than $6$ such trivial literals.
Therefore, if we have at least $11$ clauses, then we must have three such that they do not involve trivial literals in any positions.

\end{itemize}
Now it is clear that any feasible assignment of outcomes to hyperplanes in the proof of Theorem~\ref{complexity-partially-defined} can be extended to a feasible assignment of outcomes to hyperplanes for the fully defined game form by assigning outcome $*$ to $H_1^4$. Now we shall show the other direction, i.e. prove, that the assignment of outcome $*$ to $H_1^4$ is forced, i.e. there is no feasible assignment of outcomes to hyperplanes of the 4-person game form in which $H_1^4$ is assigned something else. This will in turn prove, that any feasible assignment of the fully defined 4-person game form defines a feasible assignment for the partially defined 3-person game form which is obtained by deleting all $*$ outcomes and considering the three dimensional projection by deleting the fourth coordinate (which is always 1 for all strategy profiles). Consequently, the feasible assignments of the fully defined 4-person game form correspond to satisfying assignments of the input 3-CNF $\Phi$ in the same manner as in the proof of Theorem~\ref{complexity-partially-defined}.

Let us therefore assume by contradiction that $H_1^4$ is assigned an outcome $z$ different from outcome $*$.
We have two possible cases. Either $z=c_{\ell}$ for some $1 \leq \ell \leq m$ or $z$ is one of the added outcomes participating in a gadget designed  to force values for hyperplanes incident to some $c_i$ and $c_j$ in the selector box $B$. Now let us consider three clauses $C_p, C_q, C_r$ all different from $C_{\ell}$ in the first case, or different from both $C_i$ and $C_j$ in the second case such that literals $u_p, v_q, w_r$ are nontrivial. Such clauses and literals exist due to the additional assumption on $\Phi$ stated above. This implies, that in any feasible assignment the hyperplane $H_p^1$ is forced to be assigned either $c_p$ or one of the added outcomes participating in the corresponding gadget (which are all different from $z$), and similarly for $H_q^2$ and $H_r^3$. Now observe, that the strategy profile $x$ at the intersection of $H_p^1$, $H_q^2$, and $H_r^3$ lies inside the selector box $B$ and outside its main diagonal (we assumed that $p,q,r$ are not all equal). Thus we have $g(x)=*$, and $x$ can be covered by neither $H_p^1$, nor $H_q^2$, nor $H_r^3$, nor $H_1^4$. Therefore, no feasible assignment that assigns to $H_1^4$ an outcome different from $*$ can exist.
\qed

\section{Further connections with game theory}
\label{game-theory}

In \cite{BGMP10}, the following six classes of game forms were considered:
tight ($T$), totally tight ($TT$), Nash-solvable ($NS$), dominance-solvable ($DS$),
acyclic ($AC$), and assignable ($AS$).
The classes $AS$, $T$, and $TT$ were introduced in the previous sections.
As for the remaining three classes, $AC$, $NS$, and $DS$, we refer the reader to,
for example, \cite{Mou79, Mou83, BGMP10, BCG11, Kuk11}, where these properties are
introduced both for the two-person case and for the general $n$-person case.

It was shown in \cite{BGMP10} that for the two-person case the following implications hold
\begin{equation}
\label{eq0}
AS \; \Leftarrow \; TT \; \Leftrightarrow \; AC \; \; \Rightarrow \; DS \; \Rightarrow \; NS \; \Leftrightarrow \; T.
\end{equation}
In fact, the last three implications,
$DS \; \Rightarrow \; NS \; \Leftrightarrow \; T$, were obtained long ago;
see \cite{Mou79, Mou83, GG86a} and \cite{Gur75, Gur78, Gur88}, respectively.
Furthermore, $AC \Rightarrow TT$ is obvious; see \cite{Kuk07, Kuk11} or \cite{BGMP10}.
The remaining three implications for $n=2$, that is that $TT \Rightarrow AS$, $TT \Rightarrow DS$, and $TT \Rightarrow AC$
constitute the main results of \cite{BGMP10}.
The first two of them were conjectured by Kukushkin
(\cite{Kuk07, Kuk08}, he uses the term ``separability" instead of ``assignability" \cite{Kuk11});
the last one was proven independently in \cite{Kuk07} and \cite{BGMP10}.
All three of the above implications can easily be derived from the recursive characterization of the
$TT$ two-person game forms obtained in \cite{BGMP10}.
Let us also mention that
$TT \Rightarrow NS$ easily follows from Shapley's theorem~\cite{Sha64}, although
it was not explicitly claimed there, as~\cite{Sha64} deals with zero-sum games and
their saddle points, rather than with Nash-solvable game forms.
Finally, let us note that, except for those given in~\raf{eq0},
no other implications hold between the considered six classes of game forms.
The corresponding examples to all invalid implications can be found in \cite{BGMP10}.

It was shown in \cite{BGMP10} that totally tight two-person game forms are also assignable.
Checking whether a two-person game form is TT can be done in polynomial time, as there is only a polynomial number of $2 \times 2$ subforms to check (this easy fact was observed e.g. in~\cite{BCG11} and generalized to $n$-person game forms for any $n$  in Section~\ref{complexity-section} of this paper).
This of course also applies to checking AC, since this property is in the two-person case equivalent to TT. Recognizing whether a  two-person game form is AS requires a more complicated argument, but as we have seen in Section~\ref{complexity-section} this can be done in polynomial time for $n=2$ and it is hard for $n\geq 4$.

Only a few relations are known between the properties of higher dimensional game forms. Moulin \cite{Mou83} proved that
the implications $DS \Rightarrow NS$ and $DS \Rightarrow T$ hold for every $n$, while $AC \Rightarrow NS$ is trivial.
However, tightness and Nash-solvability are no longer related, that is, both implications in $NS \; \Leftrightarrow \; T$ fail already for $n=3$.
In \cite{Gur75}, $\;\Leftarrow\;$ was disproved, while $\;\Rightarrow\;$ mistakenly claimed; then, in \cite{Gur88} it was shown that both implications fail already for $n = 3$.

%
%
%

In \cite{BCG11}, it was demonstrated that for $n\geq 3$ total tightness implies neither dominance-solvability nor acyclicity.
Both implications, $TT \Rightarrow AC$ and $TT \Rightarrow DS$, fail already for $n=3$, as shown by a single example in the last section of \cite{BCG11}.
Moreover, the same example disproves the conjecture from \cite{Kuk07} stating that a game form is acyclic if all of its subforms are Nash-solvable.
However, $TT \Rightarrow NS$, still holds for $n=3$ (this is the main result of \cite{BCG11}), while the case $n > 3$ remains open.

This last result is interesting for the following reason: while testing whether a three-person game form $g$ is Nash-solvable is computationally difficult, testing whether $g$ is (weakly) totally tight can be done in polynomial time as we have seen in Section~\ref{complexity-section}.

\section{Conclusions and plans for future research}
\label{conclusions}
This paper presents two main results for (partially defined) $n$-person game forms:
\begin{enumerate}
\item If a partially defined $n$-person game form is weakly totally tight (which is a condition verifiable in polynomial time) then the game form is assignable, and there is a polynomial time combinatorial algorithm constructing a feasible assignment.
\item It is NP-complete to decide whether a partially defined $n$-person game form is assignable for any fixed $n \geq 3$, and the same is true for fully defined $n$-person game forms for any fixed $n \geq 4$ .
\end{enumerate}
The first result of course holds also for fully defined WTT $n$-person game forms as these can be viewed as a subclass of the partially defines ones. The proofs of the second result assume that the number of outcomes $|A|$ is an input parameter. The complexity of recognition of assignable game forms for a fixed number of outcomes remains open, in particular it remains open for the interesting case $|A| = n+1$ (where $n$ is the number of players), which is the minimal number of outcomes for which non-assignable game forms exist. This case seems to be difficult, even when the last outcome appears only once. Besides of this open question, the paper introduced three conjectures concerning assignability of game forms.

\begin{itemize}
\item No-sink WTT game forms exists for any number of players.
\item For a no-sink game form, assigning the proper outcome for all hyperplanes is a feasible assignment.
\item Recognizing assignability of fully defined game forms of three players is NP-hard.
\end{itemize}

For the two-person case, the equalities and containments of \raf{eq0} summarize the relations between the corresponding six classes of game forms. Furthermore, there are no others, that is, all containments not shown by \raf{eq0} fail. These results are obtained in \cite{BGMP10} and are based on a complete recursive characterization of the $TT$ two-person game forms. The six properties of \raf{eq0} are defined for any $n$. Moreover, for  $n \geq 4$  the concepts of tightness and total tightness appear in two versions $T$, $WT$ and  $TT$, $WTT$. By definition, $T \Rightarrow WT$  and  $TT \Rightarrow WTT$.


It is a natural idea to extend the results for $n = 2$ to the general $n$-person case and complete the diagram of containments between the considered  $8$  classes of game forms. A possible plan is to obtain a recursive characterization of the $n$-person $TT$ and $WTT$ game forms and derive from it the main properties (containments) for this class in a similar manner as was done in \cite{BGMP10} for the case $n=2$.

\section*{Acknowledgement}
The first author thanks for partial support the National Science Foundation (Grant IIS-1161476). This work was partially done while the second author was a visiting professor at RUTCOR and he thanks RUTCOR for the generous support he has received. He also acknowledges the support of the Czech Science Foundation (Grant 15-15511S-P202). The research of the third author was partially funded by the Russian Academic Excellence Project '5-100'. We are also thankful to Igor Zverovich who suggested the equivalent reformulation of total tightness of Lemma \ref{constant-region}, to Vladimir Oudalov whose code was used to get the two examples in Figures \ref{f-nosink-3D-1} and \ref{f-no-sink-3D-2}, and to Nikolai Kukushkin who brought our attention to the concepts of acyclicity and assignability of game forms in his talk at RUTCOR in February 2008.

\end{document}